\documentclass[11pt, a4paper, twoside]{article}

\usepackage[left = 1in, right = 1in, top = 1.1in, bottom = 1in]{geometry}

\usepackage{amsmath, dsfont, multirow,amssymb, amsthm,mathrsfs}

\usepackage{enumerate}

\usepackage[bf]{caption}

\usepackage{graphicx,color}

\usepackage{booktabs}	

\usepackage{rotating}

\usepackage{natbib}

\usepackage[pdftex]{hyperref}
\usepackage{graphicx}
\hypersetup{
    unicode=true,
    colorlinks = true,
    linkcolor = blue,
    anchorcolor = blue,
    citecolor = blue,
    filecolor = blue,
    urlcolor = blue
}

\theoremstyle{plain}
\newtheorem{prop}{Proposition}

\theoremstyle{definition}

\newtheorem{remark}{Remark}
\newtheorem{example}{Example}

\begin{document}

 \title{\bf MCMC Bayesian Estimation in FIEGARCH Models\vspace{0.5cm}}

 \author{Taiane S. $\mbox{Prass}^{\tiny \mbox{(a)}}$\footnote{Corresponding author: taianeprass@gmail.com},
 \,  S\'ilvia R.C. $\mbox{Lopes}^{\tiny \mbox{(a)}}$ and Jorge A. $\mbox{Achcar}^{\tiny \mbox{(b)}}$\vspace{0.5cm}\\
 {\small (a)  Mathematics Institute - UFRGS, Porto Alegre, RS,  Brazil}\\
{\small \hspace{-1.5cm} (b)  Medical School  -  USP, Ribeir\~ao Preto, SP,  Brazil}\vspace{0.5cm} }

\maketitle

    \begin{abstract}
Bayesian inference for fractionally integrated exponential generalized autoregressive conditional heteroskedastic (FIEGARCH) models using Markov Chain Monte Carlo (MCMC) methods is described.   A simulation study is presented to access the performance  of the procedure,  under the presence of long-memory in the volatility.   Samples from FIEGARCH processes are obtained upon considering the generalized error distribution (GED)  for the innovation process.  Different values for the tail-thickness parameter  $\nu$  are considered covering both scenarios, innovation processes with  lighter ($\nu < 2$) and  heavier ($\nu > 2$) tails than  the  Gaussian distribution ($\nu = 2$).  A sensitivity analysis is performed by considering different prior density functions and by integrating (or not) the knowledge on the true parameter values to select  the hyperparameter values.\\


\end{abstract}

\section{Introduction}

ARCH-type (Autoregressive Conditional Heteroskedasticity)  and stochastic volatility  \citep{BEA1998} models are commonly used in financial time series modeling to represent  the dynamic evolution of volatilities.   By ARCH-type models we mean not only the ARCH model proposed by \cite{E1982}  but also  several generalizations that were lately proposed.

Among the most popular generalizations of the ARCH  model is  the  generalized ARCH (GARCH) model, introduced by \cite{B1986},  for which  the conditional variance  depends not only on the $p$ past values of the process (as in the ARCH model),  but also on the $q$ past values of the conditional variance.   Although the ARCH and GARCH models are widely used in practice,   they do not take into account  the asymmetry in the  volatility,  that is, the fact that  volatility tends to rise in response to ``bad'' news and to   fall in response to ``good'' news.   As an alternative, \cite{NE1991} introduces  the exponential GARCH (EGARCH) model.   This model not only describes the  asymmetry on the volatility,  but also have  the  advantage that the  positivity of the conditional variance is always attained since it is  defined in terms of the logarithm function.

The  fractionally integrated EGARCH (FIEGARCH) and fractionally integrated GARCH (FIGARCH) models  proposed, respectively,  by \cite{BM1996} and  \cite{BEA1996},   generalize  the EGARCH \citep{NE1991} and the GARCH \citep{B1986} models, respectively.  FIEGARCH models have not only the capability of modeling clusters of volatility (as  ARCH and GARCH models do) and capturing its asymmetry (as  the EGARCH model does) but they also take into account the characteristic of long memory in the volatility (as the FIGARCH model does).  The non-stationarity of FIGARCH models (in the weak sense) makes this class of models less attractive for practical applications.   Another drawback of the FIGARCH  models is that we must have $d \geq 0$   and the polynomial coefficients in its definition must satisfy some restrictions so the conditional variance will be positive.  FIEGARCH models do not have this problem since the variance is defined in terms of the logarithm function, moreover, they are  weak stationary whenever the long memory parameter $d$ is smaller than  0.5 \citep{LP2013}.

A complete study on the theoretical properties of FIEGARCH processes is presented in \cite{LP2013}.  The authors also conduct a simulation study to analyze the finite sample performance of the quasi-maximum  likelihood (QML) procedure on  parameter estimation.   The QML  procedure became popular  for two main reasons.  First,  the expression for the  quasi-likelihood function is simpler for the Gaussian case than when considering, for example, the Student's $t$ or the  generalized error distribution (GED).   Second,  since  the parameters  of the  distribution function are not estimated,  the dimension of  the optimization problem is reduced.     On the other hand,   the results in \cite{LP2013} indicate  that,  although the  QML presents a relatively good performance when the sample size is 2000 and the estimation improves as the sample size increases,  it does so very slowly.

In this work we propose the use of Bayesian methods  using Monte Carlo simulation techniques  on the estimation of the  FIEGARCH model parameters. This procedure is  usually considered to analyze  financial time series  assuming stochastic volatility models \citep[see, for example,][]{MY2000}, mostly because of the difficulty on applying traditional statistical techniques  due to the complexity of the likelihood function. To generate samples from  the joint posterior distribution for the parameters of interest we use MCMC (Markov Chain Monte Carlo) methods as the Gibbs Sampling algorithm  \citep[see, for example,][]{GS1990,CG1992} or the  Metropolis-Hastings algorithm \citep[see, for example,][]{SR1993,CG1995}. These samples are generated from all conditional posterior distributions for each parameter given all the other parameters and the data.

A simulation study is  conducted  to access the finite sample performance  of the procedure  proposed here, under the presence of long-memory in the volatility.   The samples from FIEGARCH processes are obtained upon considering the GED  for the innovation process.
 Taking into account that financial time series are usually characterized by heavy tailed distributions,  different values for the tail-thickness parameter  $\nu$  are considered covering both scenarios:  innovation processes with  lighter and  heavier tails than  the  Gaussian distribution.  A sensitivity analysis is performed by considering different prior density functions and by integrating (or not) the knowledge on the true parameter values to select  the hyperparameter values.

 The paper is organized as follows.  In Section \ref{section:fiegarch} a review on the definition and main properties of FIEGARCH processes is presented. Section \ref{MCMCSection} describes the parameter estimation procedure when Bayesian inference using MCMC is considered.  Section \ref{section:simulation} describes the steps used in the simulation study, such as the data generating process, the prior selection procedure and the performance measures considered.  This section also reports the simulation results. Section \ref{seciton:conclusion}  concludes  the paper.

\section{FIEGARCH Processes}\label{section:fiegarch}

Let $(1-\mathcal{B})^{-d}$ be the operator defined by its Maclaurin
series expansion, namely,
    \begin{equation}\label{binomialexp}
      (1-\mathcal{B})^{-d} = \sum_{k=0}^{\infty}\frac{\Gamma(k+d)}{\Gamma(k+1)\Gamma(d)} :=
      \sum_{k=0}^{\infty}\tau_{d,k}\,\mathcal{B}^k,
    \end{equation}
where  $\tau_{d,k} : = \frac{\Gamma(k+d)}{\Gamma(k+1)\Gamma(d)}$,  for all $k\geq1$,   $\Gamma(\cdot)$ is the gamma function and  $\mathcal{B}$ is the backward shift
operator defined by $\mathcal{B}^k (X_t) = X_{t-k}$, for all $k\in\mathds{N}$.

Assume that  $\alpha(\cdot)$ and $\beta(\cdot)$ are polynomials  of order $p$ and $q$, respectively,   defined by
    \begin{equation}
      \alpha(z) = \sum_{i=0}^{p}(-\alpha_i)z^i  \quad \mbox{ and}
      \quad \beta(z) = \sum_{j=0}^{q}(-\beta_j)z^j,\label{alphabeta}
    \end{equation}
with $\alpha_0 = \beta_0 = -1$.  If  $\alpha(\cdot)$ and $\beta(\cdot)$ have no common roots and  $\beta(z)\neq 0$ in the closed disk $\{z: |z|\leq 1\}$,  then  the function  $\lambda(\cdot)$, defined by
      \begin{equation}\label{lambdapoly}
        \lambda(z)=\frac{\alpha(z)}{\beta(z)}(1-z)^{-d} :=\sum_{k=0}^{\infty}\lambda_{d,k}z^k, \quad \mbox{ for all } |z| < 1,
    \end{equation}
is analytic in the open disc $\{z:|z|<1\}$, for any $d > 0$, and  in the closed disk $\{z:|z|\leq 1\}$, whenever $d\leq 0$.  Therefore, $\lambda(\cdot)$ is well defined and the power series representation in \eqref{lambdapoly} is unique.   More specifically,  the coefficients $\lambda_{d,k}$, for all $k \in
  \mathds{N}$, are given by \citep[see][]{LP2013}
      \begin{equation}\label{coefflambda}
        \lambda_{d,0} = 1 \quad \mbox{ and} \quad \lambda_{d,k}= -\alpha_k^* +  \sum_{i=0}^{k-1}\lambda_i \sum_{j=0}^{k-i} \beta_j^*\delta_{d,k-i-j},\
        \mbox{ for all }\ k\geq 1,
      \end{equation}
where
       \begin{equation}\label{def} \alpha_m^* := \left\{
          \begin{array}{ccc} \alpha_m, & \mbox{ if} &0
            \leq m \leq p;\vspace{.3cm}\\
            0, & \mbox{ if} & m > p;
          \end{array} \right. \quad \beta_m^* := \left\{
          \begin{array}{ccc} \beta_m, & \mbox{ if} & 0\leq
            m \leq q;\vspace{.3cm}\\
            0, & \mbox{ if} & m > q;
          \end{array} \right.
      \end{equation}
  and  $\delta_{d,j} := \tau_{-d,j}$, for all $j\in\mathds{N}$, are the coefficients obtained upon  replacing $-d$ by $d$ in \eqref{binomialexp}, that is
        \begin{equation*}
         \sum_{k=0}^{\infty}  \delta_{d,k}\mathcal{B}^{k} :=   \sum_{j=0}^{\infty} \tau_{-d,j}\mathcal{B}^{j} =  (1-\mathcal{B})^{d}.
        \end{equation*}

Let $\theta, \gamma \in \mathds{R}$ and  $\{Z_t\}_{t \in \mathds{Z}}$ be a sequence of independent and
identically distributed (i.i.d.) random variables, with zero mean and variance equal to one.
Assume that $\theta$ and $\gamma$ are not both equal to zero and define
$\{g(Z_t)\}_{t \in \mathds{Z}}$ by
      \begin{equation}\label{functiong}
        g(Z_t)=\theta Z_t +  \gamma[|Z_t|-\mathds{E}(|Z_t|)],\quad \mbox{ for all }
        t\in\mathds{Z}.
      \end{equation}
It follows that \citep[see][]{LP2013} $\{g(Z_t)\}_{t \in \mathds{Z}}$  is a strictly stationary and ergodic process. Moreover,  since $\mathds{E}(Z_0^2) <\infty$,
then   $\{g(Z_t)\}_{t\in\mathds{Z}}$ is
also weakly stationary  with mean
zero (therefore a white noise process) and variance $\sigma^2_g$ given by
    \begin{equation}\label{eq:sigmag}
      \sigma^2_g = \theta^2 + \gamma^2
      -[\gamma\mathds{E}(|Z_0|)]^2 + 2\, \theta \, \gamma \, \mathds{E}(Z_0|Z_0|).
    \end{equation}

Now, for any  $d<0.5$ and  $\omega \in \mathds{R}$, let  $\{X_t\}_{t \in \mathds{Z}}$ be the stochastic process defined by
      \begin{align}
      X_t \,& = \, \sigma_tZ_t,\\
        \ln(\sigma_t^2)&\, = \,
        \omega+\frac{\alpha(\mathcal{B})}{\beta(\mathcal{B})}(1-\mathcal{B})^{-d}g(Z_{t-1})
        \nonumber\\
       & =  \omega + \sum_{k=0}^\infty \lambda_{d,k}g(Z_{t-1-k}),\quad \mbox{ for all } t\in\mathds{Z}.\label{fieprocess}
      \end{align}
Then   $\{X_t\}_{t \in \mathds{Z}}$  is a \emph{Fractionally Integrated}  EGARCH \emph{process}, denoted by FIEGARCH$(p,d,q)$ \citep{BM1996}.\vspace{1\baselineskip}

 The properties of FIEGARCH$(p,d,q)$ processes, with $d < 0.5$, are given below  \citep[the proofs of these properties can be found in][]{LP2013}.   Henceforth  \emph{GED}$(\nu)$ denotes the generalized error distribution with tail thickness parameter $\nu$.

\begin{prop}
Let  $\{X_t\}_{t \in \mathds{Z}}$ FIEGARCH$(p,d,a)$ process. Then the following properties hold:
    \begin{enumerate}[{\rm \bf 1.}]
    \item $\{\ln(\sigma_t^2)\}_{t\in\mathds{Z}}$ is a stationary (weakly and strictly) and an ergodic process and the random  variable $\ln(\sigma_t^2)$ is almost surely finite, for all $t \in  \mathds{Z}$;

    \item  if $d \in (-1, 0.5)$ and $\alpha(z)\neq 0$, for $|z|\leq 1$, the process
      $\{\ln(\sigma_t^2)\}_{t\in\mathds{Z}}$ is   invertible;

    \item   $\{X_t\}_{t \in  \mathds{Z}}$ and $\{\sigma_t^2\}_{t\in\mathds{Z}}$ are strictly stationary
      and ergodic processes;

    \item  if  $\{Z_t\}_{t \in \mathds{Z}}$ is  a sequence of i.i.d. \emph{GED}$(\nu)$  random variables,  with  $v>1$,  zero mean and variance equal to   one,    then   $\mathds{E}(X_t^r)<\infty$ and   $\mathds{E}(\sigma_t^{2r})<\infty$, for all $t\in \mathds{Z}$ and $r>0$.
    \end{enumerate}
\end{prop}

\section{Parameter Estimation: Bayesian Inference using MCMC}\label{MCMCSection}

Let $\nu$ be the  parameter  (or vector of parameters) associated to the probability density function of $Z_0$ and denote by
    \begin{itemize}
     \item  $\boldsymbol {\eta} = (\nu, d, \theta, \gamma, \omega, \alpha_1, \cdots, \alpha_p, \beta_1, \cdots, \beta_q)' := (\eta_1, \eta_2,  \cdots, \eta_{5+p+q})'$  the vector of unknown parameters in \eqref{fieprocess};

    \item  $\boldsymbol {\eta}_{(-i)}$ the vector containing all parameters in $\boldsymbol\eta$ except  $\eta_i$, for each $i\in\{1,\cdots,5+p+q\}$;

     \item $p_{\mbox{\tiny\ensuremath{Z}}}( \cdot | \nu)$ the probability density  function of
    $Z_0$ given $\nu$;

    \item   $\mathcal{F}_{t}$ the  $\sigma$-algebra generated by $\{Z_s\}_{s\leq t}$;

    \item   $p_{\mbox{\tiny\ensuremath{X_t}}}(\cdot |\boldsymbol \eta,  \mathcal{F}_{t-1})$  the  probability density function  of $X_t$ given    $\boldsymbol {\eta}$ and $\mathcal{F}_{t-1}$, for all $t\in\mathds{Z}$.
    \end{itemize}

From \eqref{fieprocess} it is evident that,   given   $\boldsymbol \eta$,   $\sigma_t$ is a $\mathcal{F}_{t-1}$-measurable random variable. Moreover, since   $X_t = \sigma_tZ_t$  and
 $p_{\mbox{\tiny\ensuremath{Z}}}( \cdot | \nu, \mathcal{F}_{t-1}) =  p_{\mbox{\tiny\ensuremath{Z}}}( \cdot | \nu)$,  the following equality holds
    \begin{equation}
     p_{\mbox{\tiny\ensuremath{X_t}}}(x_t |\boldsymbol \eta,  \mathcal{F}_{t-1} )    =
      \frac{1}{\sigma_t}p_{\mbox{\tiny\ensuremath{Z}}}\big(x_t \sigma_t^{-1}|\nu\big),  \quad \mbox{with} \quad  \sigma_t  = \exp\bigg\{\frac{1}{2}\bigg[\omega+\sum_{k=0}^\infty\lambda_{d,k}g(z_{t-1-k})\bigg]\bigg\}, \label{eq2}
    \end{equation}
for all $x_t\in\mathds{R}$ and $t\in\mathds{Z}$.
Furthermore, from \eqref{eq2},  the  conditional probability of $\boldsymbol{X} := (X_1, \cdots, X_n)'$  given  $\boldsymbol \eta$ and $ \mathcal{F}_{0}$  can be written as
    \begin{align}
    p_{\mbox{\tiny\ensuremath{\boldsymbol{X}}}}\big (x_1, \cdots, x_n | \boldsymbol{\eta}, \mathcal{F}_{0} \big) &=  p_{\mbox{\tiny\ensuremath{X_n}}}(x_n |\boldsymbol \eta, x_{n-1}, \cdots, x_1, \mathcal{F}_{0} )\times\cdots\times p_{\mbox{\tiny\ensuremath{X_1}}}(x_1 |\boldsymbol \eta,  \mathcal{F}_{0} ) \nonumber\\
    & = \prod_{t=1}^n\frac{1}{\sigma_t}p_{\mbox{\tiny\ensuremath{Z}}}\big(x_t \sigma_t^{-1}|\nu\big).\label{conditionalX}
      \end{align}

Given  any  $I_0 \in \mathcal{F}_{0}$, select  a prior  conditional  density function $p_{\mbox{\tiny\ensuremath{I_0}}}(\cdot|\boldsymbol \eta)$ for $I_0$ given $\boldsymbol \eta$. Also, select  a prior\footnote{In fact, the priors $\pi_i(\cdot)$ are not necessarily probability density functions. For instance, $\pi(x)=1$ and $\pi(x) = 1/x$,  are  examples of  improper  priors (i.e., they do not integrate to 1) used in practice.}   density function $\pi_i(\cdot)$  for  $\eta_i$  and  a prior conditional  probability density  function $p_{\mbox{\tiny\ensuremath{(-i)}}}(\cdot| \eta_i)$ for $\boldsymbol{\eta}_{(-i)}$ given $\eta_i$,  for each $i\in \{1,\cdots, 5+p+q\}$.

Observe that,  by applying the Bayes'  rule, the  conditional probability density function of $\eta_i$ given $\boldsymbol{X}$,  $\boldsymbol{\eta}_{(-i)}$ and any  $I_0$, can be written as
    \begin{equation}
    p\big(\eta_i | \boldsymbol{X}, \boldsymbol{\eta}_{(-i)}, I_{0}\big)  =
    \frac{ p_{\mbox{\tiny\ensuremath{\boldsymbol{X}}}} \big(\boldsymbol{X} | \boldsymbol{\eta}, I_{0} \big)
    \times  p_{\mbox{\tiny\ensuremath{I_0}}}\big(I_0|\boldsymbol \eta \big) \times p_{\mbox{\tiny\ensuremath{(-i)}}}(\boldsymbol{\eta}_{(-i)}| \eta_i) \times\pi_i(\eta_i)}{p_{\mbox{\tiny\ensuremath{(-i)}}}\big(\boldsymbol{X},\boldsymbol{\eta}_{(-i)}, I_{0} \big)}, \label{peta_i}
    \end{equation}
for each $i\in\{1,\cdots, 5+p+q\}$, where $p_{\mbox{\tiny\ensuremath{\boldsymbol{X}}}}\big (\cdot| \boldsymbol{\eta}, \mathcal{F}_{0} \big)$ is  given in  \eqref{conditionalX} and $p_{\mbox{\tiny\ensuremath{(-i)}}}\big(\cdot,\cdot,\cdot\big )$ is the joint probability density function of $\boldsymbol{X}$, $\boldsymbol{\eta}_{(-i)}$ and  $I_{0}$, which does not depend on $\eta_i$.

The parameter estimation is then carried out by  using the MCMC method as described below.

\subsection{Gibbs Sampling with Metropolis Steps}

Gibbs sampling \citep{GG1984,GS1990} is a popular MCMC algorithm for obtaining a sequence of random samples from multivariate probability distribution when direct sampling is difficult.  The algorithm assumes that   the conditional distribution of each random variable  is known and it is easy to sample from it.  The steps of the sampling procedure are the following.

    \begin{enumerate}[{\bf Step 1.}]
    \item Set an arbitrary initial value for the vector of parameters $\boldsymbol \eta$, namely, \\
    $\boldsymbol \eta^{(0)} = (\eta_1^{(0)}, \cdots, \eta_{5+p+q}^{(0)})'$.  Let $m = 0$;

    \item Given the sample $\boldsymbol \eta^{(m)} = (\eta_1^{(m)},\cdots, \eta_{5+p+q}^{(m)})'$,
        \begin{enumerate}[--]
        \item  generate \, $\eta_1^{(m+1)}$ \, from \,\, $p\big(\eta_1 | \boldsymbol{X}, \eta_2^{(m)}, \eta_3^{(m)}, \cdots, \eta_{5+p+q}^{(m)}, I_{0}\big)$;
        \item  generate \, $\eta_2^{(m+1)}$ \, from \,\, $p\big(\eta_2| \boldsymbol{X}, \eta_1^{(m+1)}, \eta_3^{(m)}, \cdots, \eta_{5+p+q}^{(m)}, I_{0}\big)$;\\
         \phantom{x}\hspace{2cm}$\vdots$ \hspace{2cm}$\vdots$ \hspace{2cm}$\vdots$\\
        \item  generate \, $\eta_{5+p+q}^{(m+1)}$ \, from \,\, $p\big(\eta_{5+p+q} | \boldsymbol{X}, \eta_1^{(m+1)}, \cdots \eta_{4+p+q}^{(m+1)}, I_{0}\big)$;
        \end{enumerate}

    \item Once  the vector  $\boldsymbol \eta^{(m+1)} = (\eta_1^{(m+1)}, \cdots, \eta_{5+p+q}^{(m+1)})'$  is obtained, return to step 2,  with $m = m+1$, until $m = N$,  where $N$ is the desired sample size.
    \end{enumerate}

When it is not possible to sample directly  from $p\big(\eta_i | \boldsymbol{X}, \boldsymbol{\eta}_{(-i)}, I_{0}\big)$, for one or more $i\in\{1,\cdots, 5+p+q\}$, an alternative option is to consider a combination of Gibbs sampler and Metropolis-Hastings \citep{M1953, H1970}  algorithms. This method is usually referred to as Gibbs sampler with Metropolis steps.  In this case,  to draw the random variate $\eta_i$,  one shall follow the same steps 1-3 just described. However,   instead of sampling  directly from  $p\big(\eta_i | \boldsymbol{X}, \boldsymbol{\eta}_{(-i)}, \mathcal{F}_{0}\big)$, one shall consider the Metropolis-Hastings algorithm with  $p\big(\eta_i | \boldsymbol{X}, \boldsymbol{\eta}_{(-i)}, \mathcal{F}_{0}\big)$ as the invariant (target) distribution.

Metropolis-Hastings algorithm is easy to implement since it  does not require knowing the  normalization constant $p_{\mbox{\tiny\ensuremath{(-i)}}}\big(\boldsymbol{X},\boldsymbol{\eta}_{(-i)}, I_{0} \big)$, defined  in \eqref{peta_i}.     For simplicity of notation,   in what follows  $p_*(\cdot)$ shall denote any one of  the non-normalized probability density function which corresponds  to  $p\big(\eta_i | \boldsymbol{X}, \boldsymbol{\eta}_{(-i)}, \mathcal{F}_{0}\big)$, for  $i\in\{1,\cdots, 5+p+q\}$.   The Metropolis-Hastings sampling procedure consists if the following steps.

    \begin{enumerate}[{\bf Step 1.}]
    \item Select a transition kernel\footnotemark (also called proposal distribution) $q(\cdot| \cdot)$  for which the sampling procedure is known.

    \item  Set an arbitrary initial value $y_0$ for the chain. Let $m = 0$.

    \item  Generate a draw $\xi$ from $q(\cdot|y_m)$.

    \item Calculate $\alpha(y_m, \xi) =  \displaystyle\min\bigg\{1, \frac{p_*(\xi)q(y_m|\xi)}{p_*(y_m)q( \xi | y_m)}\bigg\}.$

    \item Draw $u \sim  \mathcal{U} [0, 1]$.

    \item  Define
    $y_{m+1} =\left\{
    \begin{array}{cl}
    \xi,  & \mbox{if} \,\, u < \alpha(y_m , \xi);\\
    y_m, & \mbox{otherwhise}.
    \end{array}
    \right.$

    \item If $m + 1 <  N$ (where $N$  is the desired sample size),  let $m = m+1$  and go to  Step 3.
    \end{enumerate}

\footnotetext{A transition kernel is a  function  $q(x|y)$ which is a probability measure with respect to $x$,  so $\int q(x|y) dx = 1$.}

    \begin{remark}$\phantom{x}$\vspace{-1\baselineskip}

    \begin{enumerate}[{\bf 1.}]

    \item   When considering Gibbs sampler with Metropolis steps only one iteration of Metropolis-Hastings algorithm is performed for each Gibbs sampler iteration.

    \item  In both cases,  Gibbs sampler  and Metropolis-Hastings algorithm, it is advised to discard the first $B$ (for some $B < N$) observations (that is, the burn-in sample) to assure the chain convergence.

    \item   The sample obtained from the algorithm described above  is not independent.   An alternative  is  to run parallel chains instead. Another common     strategy to reduce sample autocorrelations is thinning  the Markov chain, that is, to
    keep only every $k$-th simulated draw from each sequence.  There is some controversy surrounding
    the question of whether or not it is better to run one long chain or several shorter ones \citep{GR1992, G1992}. Also,  \cite{MB1994}  show that one always get more precise posterior estimates if the entire Markov chain is used instead of the thinned one.

    \end{enumerate}
    \end{remark}

\section{Simulation Study}\label{section:simulation}

 This simulation study considers  FIEGARCH$(0,d,0)$ processes. Under this scenario, the vector of unknown parameters is  $\boldsymbol {\eta} =  (\nu, d, \theta, \gamma, \omega)':= (\eta_1, \cdots, \eta_5)' $.
The Bayesian inference approach, using MCMC to obtain posterior density functions, is used to estimate the parameters of the model.

\subsection{Data Generating Process}
The samples from  FIEGARCH$(0,d,0)$ processes are obtained by   setting  the following.

    \begin{itemize}
    \item $Z_0 \sim \mbox{GED}(\nu)$, with zero mean and variance equal to one. Thus,

    \begin{equation}
    p_{\mbox{\tiny\ensuremath{Z}}}( z | \nu)=
    \frac{\nu\exp\big\{-\frac{1}{2}|z\lambda_\nu^{-1}|^\nu\big\}}{\lambda_\nu 2^{1+1/\nu}\Gamma(1/\nu)}, \quad
    \lambda_\nu =\bigg[2^{-2/\nu}\frac{\Gamma(1/\nu)}{\Gamma(3/\nu)}\bigg]^{1/2}\!\!\!\!\!\!\!,
    \quad \quad \mbox{for all } z\in\mathds{R};\nonumber
    \end{equation}

    \item     $d\in\{0.10, \, 0.25, \, 0.35, \, 0.45\}$ and $\nu \in \{1.1, 1.5, 1.9, 2.5, 5\}$;

    \item for all models,  $\omega = -5.40$, $\theta = -0.15$ and $\gamma =  0.24$.  These values are close to the ones already observed in practical applications \citep[see, for instance,][]{NE1991, BM1996, RV2008, LP2013}.

    \item   the infinite sum in  \eqref{fieprocess} is truncated at $m^* = 50,000$.
    \end{itemize}

For each combination of $d$ and $\nu$,   a sample  $\{z_t\}_{t=-m^*}^n$, of size  $m^*+n+1$, is drawn from the GED($\nu$) distribution and then the
 sample $\{x_t\}_{t=1}^n$, from the FIEGARCH$(0,d,0)$ process, is obtained through the relation
    \[
    \quad \ln(\sigma_t^2) = \omega +
    \sum_{k=0}^{m^*}\lambda_{d,k}g(z_{t-1-k})  \quad \mbox{and } \quad x_t = \sigma_tz_t,  \quad \mbox{for all } t=1,\cdots, n.
    \]

\subsection{Parameter Estimation Settings}\label{section:priors}

The samples from the posterior distributions are obtained by considering  the Gibbs sampler algorithm with Metropolis steps as  described in Section \ref{MCMCSection}.  The  transition kernel  $q(\cdot|\cdot)$ considered in the   Metropolis-Hastings algorithm  is the function defined as
    \[
     q(x|y) = f(x; y, \sigma, a, b),
     \]
 where $f(\cdot; \cdot, \cdot, \cdot, \cdot)$ is the truncated normal density function, defined as
    \[
    f(x; \mu, \sigma, a, b) =
    \left\{
    \begin{array}{cl}
    \displaystyle\frac{1}{\sigma}\frac{\phi\big(\frac{x-\mu}{\sigma}\big)}{\Phi\big(\frac{b-\mu}{\sigma}\big) - \Phi\big(\frac{a-\mu}{\sigma}\big) }, & \mbox{if } a \leq x \leq b, \vspace{0.3cm}\\
    0, & \mbox{otherwise},
    \end{array}
    \right.
    \]
where $\phi(\cdot)$ and $\Phi(\cdot)$ are, respectively,  the probability density  and  cumulative distribution functions of the standard normal distribution;
 $a,b \in \mathds{R}$ are, respectively, the lower and upper limits of the distribution's support;  $\mu$  and $\sigma$ denote, respectively,   the distribution's (non-truncated version).

To select a reasonable  $\boldsymbol \eta^{(0)}$,  $p_{\mbox{\tiny\ensuremath{\boldsymbol{X}}}}\big (\boldsymbol{X} | \boldsymbol{\eta}, \mathcal{F}_{0} \big)$ is calculated  for different combinations of $\nu, d, \theta, \gamma$ and $\omega$. Then   $\boldsymbol \eta^{(0)}$ is defined as the vector $\boldsymbol \eta = (\nu, d, \theta, \gamma, \omega)'$ with higher likelihood function value.  To eliminate any dependence on the initial $\boldsymbol \eta^{(0)}$  a burn-in of size 1000 is considered.

A sample obtained by the method being described will probably present significative correlation\footnote{In fact,   for the parameter $d$, this correlation could  only be removed when  the  thinning parameter $\mathfrak{t}$ was set  to 200.}. However,  due to  the ergodicity property of the Markov chain, the estimation of the mean is  not  affect by the correlation in the sample.  Therefore,  to avoid  unnecessary computational work, which ultimately would not lead to improvement in terms of parameter estimation,   thinning is not implemented.    Nevertheless, an example showing the influence of using the entire chain,  the  thinned chain or only the first 1000 observations of the entire chain  (after burn-in)  is provided in the following.

 \begin{example}
Let  $\boldsymbol {\eta} =  (\nu, d, \theta, \gamma, \omega)':= (\eta_1, \cdots, \eta_5)' $ and assume that
\begin{equation}
\pi_i(\eta_i) = \left\{\begin{array}{cl}
c_i, & \mbox{if } \eta_i \in  I_i;\\
0, & \mbox{otherwise};
\end{array}
\right.\quad \mbox{for each } \quad i\in\{1,\cdots,5\},\label{eq:example}
\end{equation}
with $c_1 = 1$, $c_2 = c_3 = c_4 =  2$, $c_5 = 1/30$,   $I_1 = (0, \infty)$, $I_2 = [0, 0.5]$, $I_3[-0.5, 0]$, $I_4 = I_2$ and $I_5 = [-15, 15]$.

In the sequel, $\{\eta_i^{(k)}\}_{k=1}^{\mathfrak{n}}$ denotes the chain  of size $\mathfrak{n}$ obtained  from the posterior distribution of  $\eta_i$,   upon considering the prior $\pi_i(\eta_i)$ defined in \eqref{eq:example},  for each  $i\in \{1,\cdots, 5\}$.  Also,   $\mathfrak{b}$,  $\mathfrak{t}$ and  $N$ denote, respectively,  the burn-in size,  the thinning parameter and  the sample size of the thinned chain\footnote{Observe that,  by setting  $\mathfrak{b} = 1000$ and $\mathfrak{t} = 200$,  then a thinned chain of size  $N=1000$ can only be obtained from $\{\eta_i^{(k)}\}_{k=1}^{\mathfrak{n}}$ when $\mathfrak{n} \geq \mathfrak{b} + 1 + \mathfrak{t}(N-1) =$ 200,801.} obtained from $\{\eta_i^{(k)}\}_{k=1}^{\mathfrak{n}}$, for any $i\in \{1,\cdots, 5\}$.

Figure  \ref{figure:sample}  presents the graph of $\{\eta_i^{(k)}\}_{k=1}^{\mathfrak{n}}$, for each  $i\in \{1,\cdots, 5\}$,  with  $\mathfrak{n} = $ 200,801. Figure \ref{figure:sample} also shows  the thinned chain of size $N=1000$ obtained  by considering  $\mathfrak{b} = 1000$ and  $\mathfrak{t} = 200$. Furthermore, Figure \ref{figure:sample} gives  the sample of size 1000,  obtained  from $\{\eta_i^{(k)}\}_{k=1}^{\mathfrak{n}}$  by considering  a burn-in  equal to  1000 and no thinning, for each  $i\in \{1,\cdots, 5\}$.    The true parameter values of the FIEGARCH$(0,d,0)$ model corresponding to  these graphs are $\nu_0 = 1.9 $, $d_0 = 0.25 $, $\theta_0 =  -0.15 $,  $\gamma_0 =  0.24 $ and $\omega_0 =  -5.4 $.  Figure  \ref{figure:densities} gives the histogram and kernel density  functions corresponding to each sample in Figure  \ref{figure:sample}.  The graphs of the prior  $\pi_i(\eta_i)$  defined in \eqref{eq:example},  for  $i\in \{1,\cdots, 5\}$,  are represented in Figure \ref{figure:densities} by the dashed lines.   For a better visualization of the posterior distributions, in Figure \ref{figure:densities},  the range for the $x$-axis was restricted to the intervals $[-1.5, 2.5]$, $[-0.5,0]$, $[0, 0.5]$ and  $[-5.6, -5.1]$, respectively, for the parameters $\nu, \theta, \gamma$ and $\omega$.

\begin{figure}[!ht]
\centering
\setlength{\tabcolsep}{2pt}
\begin{tabular}{ccccc}
$\nu$ & $d$ & $\theta$ & $\gamma$ & $\omega$\\
\includegraphics[width = 0.19\textwidth]{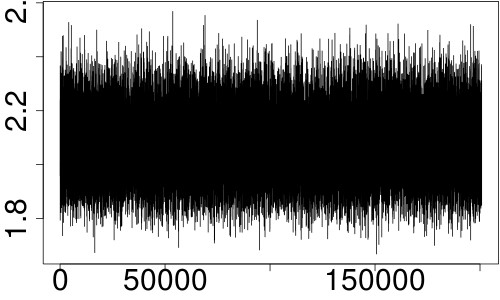} &
\includegraphics[width = 0.19\textwidth]{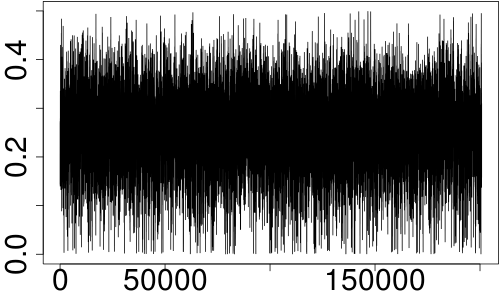} &
\includegraphics[width = 0.19\textwidth]{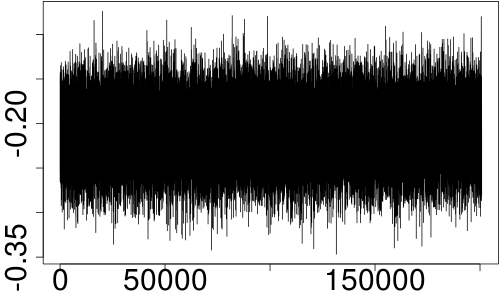} &
\includegraphics[width = 0.19\textwidth]{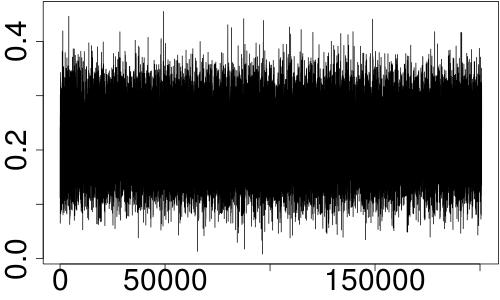} &
\includegraphics[width = 0.19\textwidth]{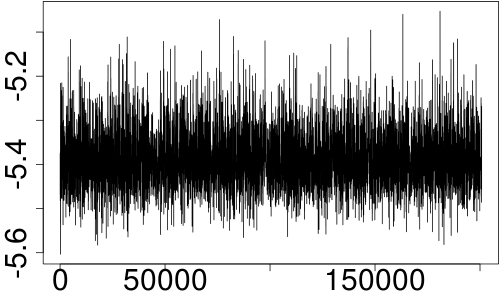} \\

\includegraphics[width = 0.19\textwidth]{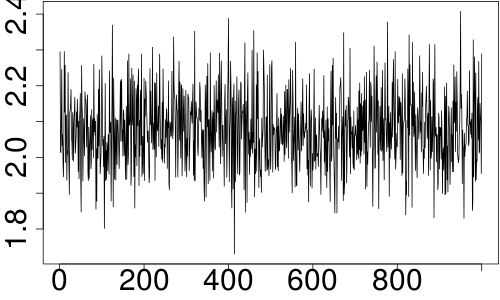} &
\includegraphics[width = 0.19\textwidth]{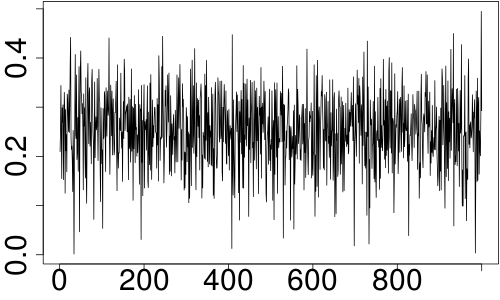} &
\includegraphics[width = 0.19\textwidth]{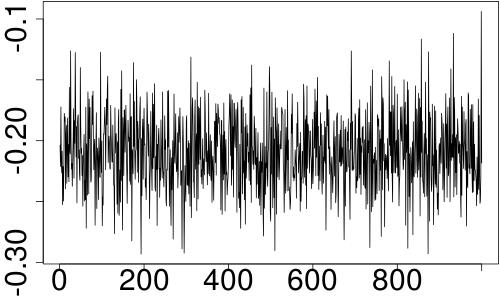} &
\includegraphics[width = 0.19\textwidth]{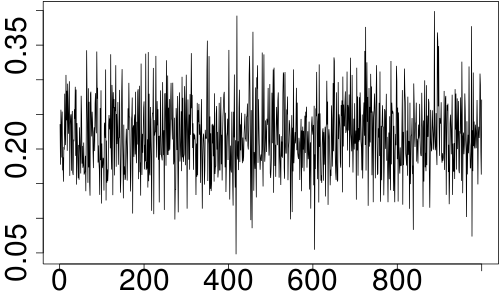} &
\includegraphics[width = 0.19\textwidth]{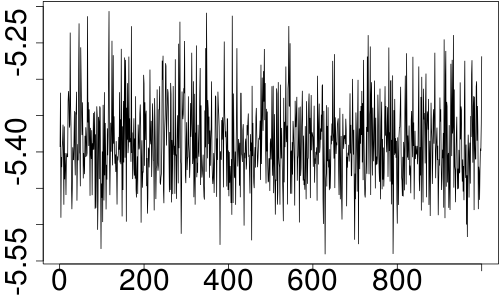} \\

\includegraphics[width = 0.19\textwidth]{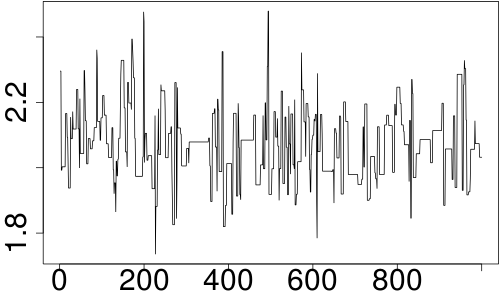} &
\includegraphics[width = 0.19\textwidth]{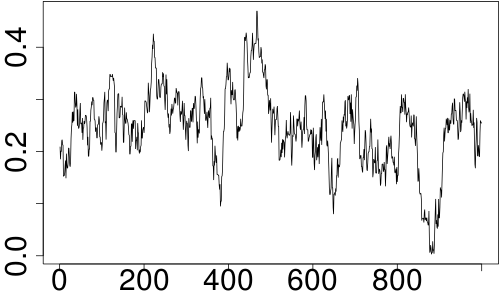} &
\includegraphics[width = 0.19\textwidth]{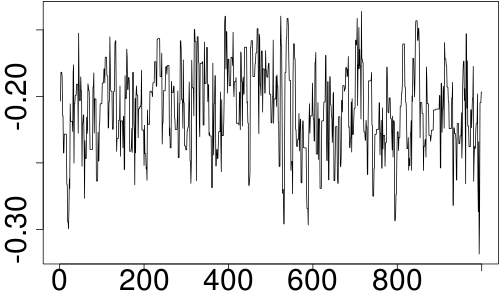} &
\includegraphics[width = 0.19\textwidth]{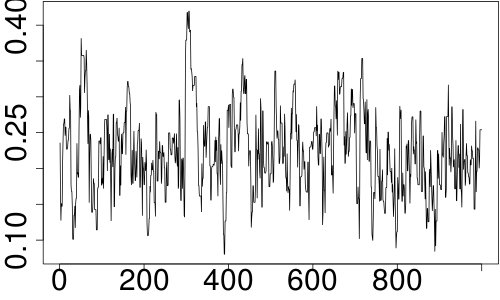} &
\includegraphics[width = 0.19\textwidth]{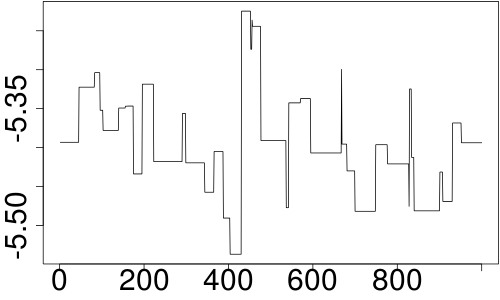} \\

\end{tabular}
\caption{Original chain with sample size 200801 (top row).
        Thinned chain with sample size 1000 and thinning parameter equal to 200 (middle row).
        Unthinned chain with sample size 1000 (bottom row).
        For  the middle and bottom rows the burn-in size is equal to 1000.
        The true parameter values of the FIEGARCH$(0,d,0)$ model corresponding to  these graphs are $\nu_0 = 1.9 $, $d_0 = 0.25 $, $\theta_0 =  -0.15 $,  $\gamma_0 =  0.24 $ and $\omega_0 =  -5.4 $.}\label{figure:sample}
\end{figure}

\begin{figure}[!ht]
\centering
\setlength{\tabcolsep}{2pt}
\begin{tabular}{ccccc}
$\nu$ & $d$ & $\theta$ & $\gamma$ & $\omega$\\
\includegraphics[width = 0.19\textwidth]{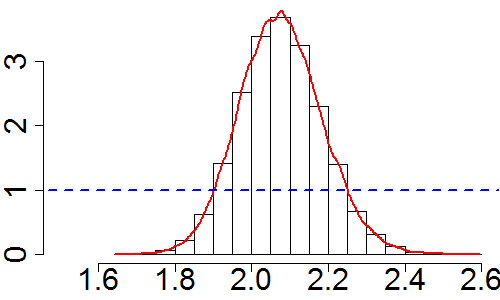} &
\includegraphics[width = 0.19\textwidth]{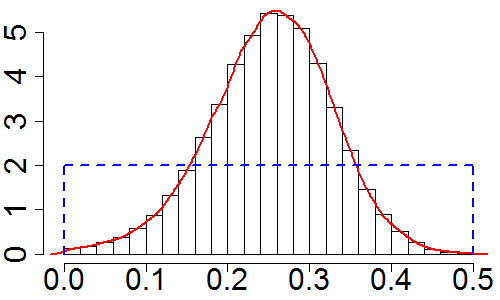} &
\includegraphics[width = 0.19\textwidth]{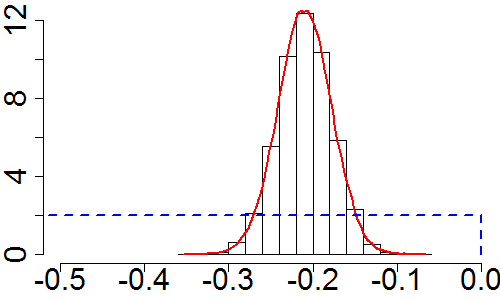} &
\includegraphics[width = 0.19\textwidth]{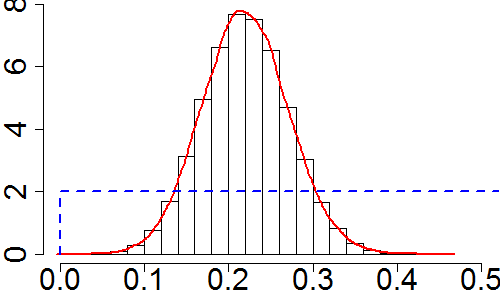} &
\includegraphics[width = 0.19\textwidth]{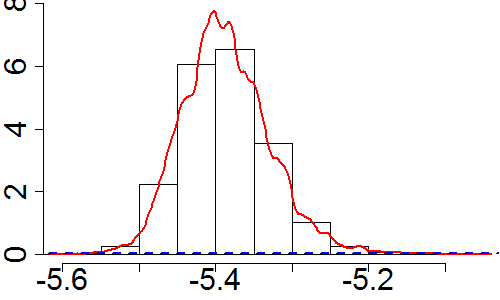} \\

\includegraphics[width = 0.19\textwidth]{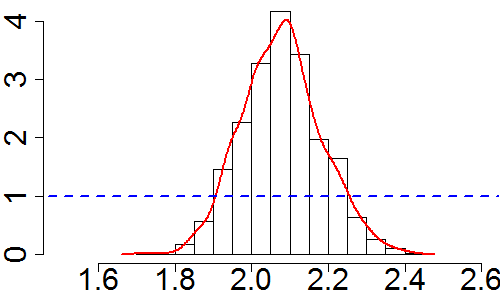} &
\includegraphics[width = 0.19\textwidth]{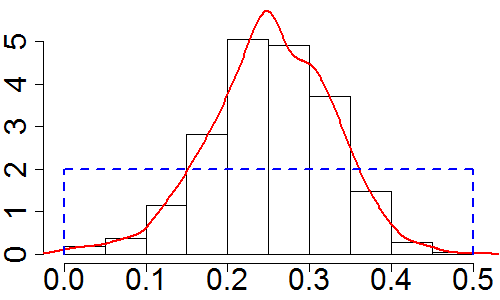} &
\includegraphics[width = 0.19\textwidth]{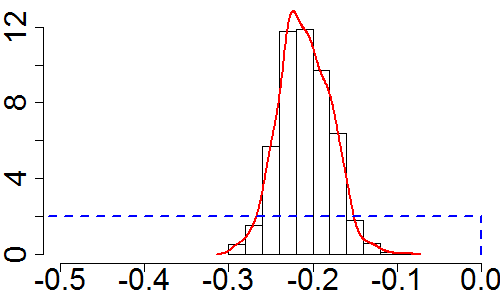} &
\includegraphics[width = 0.19\textwidth]{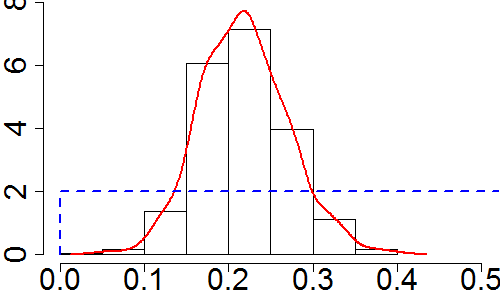} &
\includegraphics[width = 0.19\textwidth]{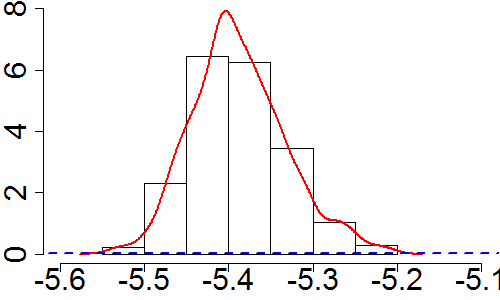} \\

\includegraphics[width = 0.19\textwidth]{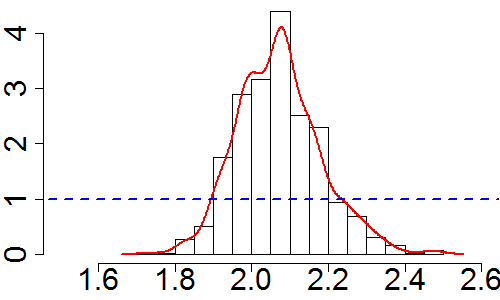} &
\includegraphics[width = 0.19\textwidth]{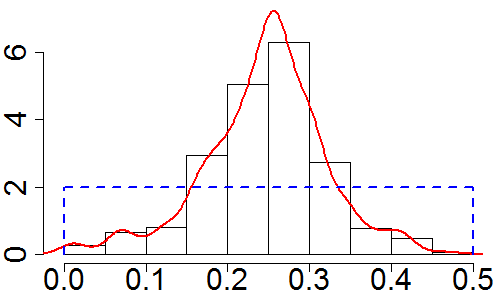} &
\includegraphics[width = 0.19\textwidth]{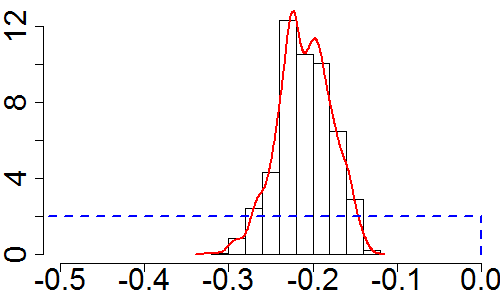} &
\includegraphics[width = 0.19\textwidth]{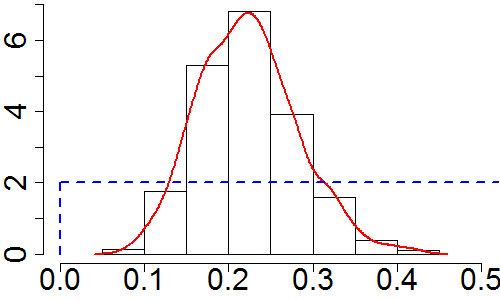} &
\includegraphics[width = 0.19\textwidth]{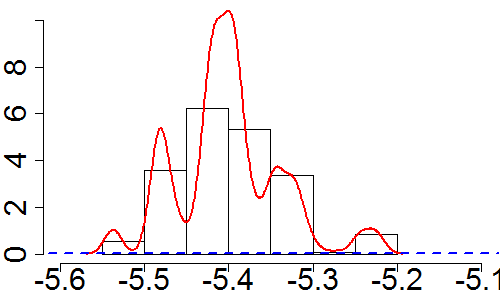} \\

\end{tabular}
\caption{Histogram and kernel density functions for the original chain with sample size 200801 (top row);
        the thinned chain with sample size 1000 and thinning parameter equal to 200 (middle row) and
        the unthinned chain with sample size 1000 (bottom row).
        For  the middle and bottom rows the burn-in size is equal to 1000.
        The dashed lines correspond to the graphs of the priors  $\pi_i(\eta_i)$  defined in \eqref{eq:example},  for  $i\in \{1,\cdots, 5\}$.   The range for the $x$-axis was restricted to the intervals $[-1.5, 2.5]$, $[-0.5,0]$, $[0, 0.5]$ and  $[-5.6, -5.1]$, respectively, for the parameters $\nu, \theta, \gamma$ and $\omega$.  The true parameter values of the FIEGARCH$(0,d,0)$ model corresponding to  these graphs are $\nu_0 = 1.9 $, $d_0 = 0.25 $, $\theta_0 =  -0.15 $,  $\gamma_0 =  0.24 $ and $\omega_0 =  -5.4 $.}\label{figure:densities}
\end{figure}

As shown in Figure \ref{figure:sample} (see also Table \ref{table:example}),  the mean of the posterior distribution does not change significantly when the entire sample or the  thinned chain is considered instead of the unthinned one.   On the other hand, Figure \ref{figure:densities} reinforces the idea that the entire chain gives better estimates  for the density function (notice that the curves in the graphs are smoother).   Although the thinned chain is not as efficient as the entire chain, it still provides better estimates for the density function than the unthinned one.

Table \ref{table:example}  presents the  summary statistics  for the samples obtained from  the posterior distribution  for  each parameter of the  FIEGARCH$(0,d,0)$ model.  This table considers  the entire, thinned and unthinned chains.   The statistics reported in this table are the sample  mean ($\bar \eta_i$), the sample standard deviation ($\mbox{sd}_{\eta_i}$) and the 95\% credibility interval $CI_{0.95}(\eta_i)$ for the parameter $\eta_i$  in $\boldsymbol {\eta} =  (\nu, d, \theta, \gamma, \omega)' :=  (\eta_1, \cdots, \eta_5)' $, for each  $i\in \{1,\cdots, 5\}$.    The true parameter values considered for this illustration are $\nu_0 \in \{1.1, 1.5, 2.5, 5.0\}$, $d_0 = 0.25$, $\theta_0 = -0.15$, $\gamma_0 = 0.24$ and $\omega_0 = -5.4$.

\begin{table}[!ht]
\centering
\setlength{\tabcolsep}{5pt}
{\footnotesize
\caption{Summary for the entire,  thinned (thinning parameter $\mathfrak{t} = 200$) and unthinned sample  from posterior distributions considering all prior uniforms:  mean $\bar \eta_i$, standard deviation $\mbox{sd}_{\eta_i}$ and the 95\% credibility interval $CI_{0.95}(\eta_i)$ for the parameter $\eta_i$  in $\boldsymbol {\eta} =  (\nu, d, \theta, \gamma, \omega)' :=  (\eta_1, \cdots, \eta_5)' $, for each  $i\in \{1,\cdots, 5\}$.  The true parameter values considered in this simulation are  $\nu_0 \in \{1.1, 1.5, 2.5, 5.0\}$, $d_0 = 0.25$, $\theta_0 = -0.15$, $\gamma_0 = -0.24$ and $\omega_0 = -5.4$. For the thinned and unthinned samples, the burn-in size is $\mathfrak{b} = 1000$.}\label{table:example}
\begin{tabular}{c|ccccccccccc}
\toprule
\multicolumn{1}{c}{\multirow{2}{*}{Chain}}&
\multicolumn{1}{c}{\multirow{2}{*}{$\nu_0$}}&&
\multicolumn{1}{c}{$\bar{\nu}$ $(\mbox{sd}_{\nu})$} &&
\multicolumn{1}{c}{$\bar{d}$ $(\mbox{sd}_{d})$} &&
\multicolumn{1}{c}{$\bar{\theta}$ $(\mbox{sd}_{\theta})$} &&
\multicolumn{1}{c}{$\bar{\gamma}$ $(\mbox{sd}_{\gamma})$} &&
\multicolumn{1}{c}{$\bar{\omega}$ $(\mbox{sd}_{\omega})$} \\
\multicolumn{1}{c}{}&\multicolumn{1}{c}{}&&
\multicolumn{1}{c}{$CI_{0.95}(\nu)$} &&
\multicolumn{1}{c}{$CI_{0.95}(d)$} &&
\multicolumn{1}{c}{$CI_{0.95}(\theta)$} &&
\multicolumn{1}{c}{$CI_{0.95}(\gamma)$} &&
\multicolumn{1}{c}{$CI_{0.95}(\omega)$} \\

\cmidrule{1-2} \cmidrule{4-4} \cmidrule{6-6} \cmidrule{8-8}\cmidrule{10-10}\cmidrule{12-12}
\multirow{10}{*}{\begin{sideways}Entire\end{sideways}}
 & \multirow{2}{*}{1.1} && 1.095 (0.047) &&0.263 (0.112) &&-0.087 (0.038) &&0.233 (0.062) &&-5.469 (0.078) \\
&&& [1.006; 1.188] && [0.035; 0.464] && [-0.164; -0.016] && [0.114; 0.358] && [-5.612; -5.299] \\
\cmidrule{2-2} \cmidrule{4-4} \cmidrule{6-6} \cmidrule{8-8}\cmidrule{10-10}\cmidrule{12-12}
 & \multirow{2}{*}{1.5} && 1.478 (0.067) &&0.220 (0.081) &&-0.184 (0.037) &&0.240 (0.057) &&-5.408 (0.058) \\
&&& [1.351; 1.612] && [0.056; 0.375] && [-0.257; -0.111] && [0.130; 0.355] && [-5.520; -5.291] \\
\cmidrule{2-2} \cmidrule{4-4} \cmidrule{6-6} \cmidrule{8-8}\cmidrule{10-10}\cmidrule{12-12}
 & \multirow{2}{*}{1.9} && 2.077 (0.107) &&0.252 (0.075) && -0.209 (0.032) &&0.220 (0.051) &&-5.386 (0.058) \\
&&& [1.874; 2.297] && [0.092; 0.392] && [-0.272; -0.148] && [0.122; 0.322] && [-5.487; -5.261] \\
\cmidrule{2-2} \cmidrule{4-4} \cmidrule{6-6} \cmidrule{8-8}\cmidrule{10-10}\cmidrule{12-12}
 & \multirow{2}{*}{2.5} && 2.727 (0.153) &&0.298 (0.056) &&-0.203 (0.025) &&0.205 (0.045) &&-5.361 (0.053) \\
&&& [2.441; 3.040] && [0.184; 0.405] && [-0.253; -0.153] && [0.118; 0.296] && [-5.463; -5.252] \\
\cmidrule{2-2} \cmidrule{4-4} \cmidrule{6-6} \cmidrule{8-8}\cmidrule{10-10}\cmidrule{12-12}
 & \multirow{2}{*}{5.0} && 5.227 (0.366) &&0.220 (0.051) &&-0.173 (0.019) &&0.294 (0.036) &&-5.303 (0.039) \\
&&& [4.548; 5.978] && [0.115; 0.317] && [-0.211; -0.135] && [0.224; 0.366] && [-5.384; -5.230] \\
\midrule

\multirow{10}{*}{\begin{sideways}Thinned\end{sideways}}
 & \multirow{2}{*}{1.1} && 1.095 (0.047) &&0.264 (0.110) &&-0.087 (0.038) &&0.233 (0.062) &&-5.469 (0.080) \\
&&& [1.006; 1.193] && [0.038; 0.462] && [-0.161; -0.015] && [0.116; 0.362] && [-5.609; -5.288] \\
\cmidrule{2-2} \cmidrule{4-4} \cmidrule{6-6} \cmidrule{8-8}\cmidrule{10-10}\cmidrule{12-12}
 & \multirow{2}{*}{1.5} && 1.480 (0.067) &&0.221 (0.077) &&-0.184 (0.037) &&0.240 (0.057) &&-5.405 (0.060) \\
&&& [1.353; 1.613] && [0.070; 0.367] && [-0.257; -0.104] && [0.129; 0.352] && [-5.524; -5.291] \\
\cmidrule{2-2} \cmidrule{4-4} \cmidrule{6-6} \cmidrule{8-8}\cmidrule{10-10}\cmidrule{12-12}
 & \multirow{2}{*}{1.9} && 2.079 (0.103) &&0.253 (0.075) &&-0.210 (0.031) &&0.218 (0.052) &&-5.388 (0.056) \\
&&& [1.888; 2.292] && [0.098; 0.387] && [-0.269; -0.152] && [0.120; 0.327] && [-5.486; -5.261] \\
\cmidrule{2-2} \cmidrule{4-4} \cmidrule{6-6} \cmidrule{8-8}\cmidrule{10-10}\cmidrule{12-12}
 & \multirow{2}{*}{2.5} && 2.718 (0.154) &&0.298 (0.058) &&-0.202 (0.026) &&0.205 (0.044) &&-5.361 (0.053) \\
&&& [2.427; 3.064] && [0.181; 0.410] && [-0.255; -0.153] && [0.121; 0.291] && [-5.466; -5.254] \\
\cmidrule{2-2} \cmidrule{4-4} \cmidrule{6-6} \cmidrule{8-8}\cmidrule{10-10}\cmidrule{12-12}
 & \multirow{2}{*}{5.0} && 5.224 (0.363) &&0.218 (0.052) &&-0.173 (0.019) &&0.295 (0.037) &&-5.302 (0.039) \\
&&& [4.534; 5.991] && [0.112; 0.316] && [-0.211; -0.137] && [0.224; 0.370] && [-5.391; -5.230] \\
\midrule

\multirow{10}{*}{\begin{sideways}Unthinned\end{sideways}}
 & \multirow{2}{*}{1.1} && 1.108 (0.039) &&0.265 (0.129) &&-0.089 (0.038) &&0.230 (0.060) &&-5.476 (0.052) \\
&&& [1.028; 1.198] && [0.016; 0.467] && [-0.176; -0.019] && [0.108; 0.345] && [-5.553; -5.366] \\
\cmidrule{2-2} \cmidrule{4-4} \cmidrule{6-6} \cmidrule{8-8}\cmidrule{10-10}\cmidrule{12-12}
 & \multirow{2}{*}{1.5} && 1.474 (0.067) &&0.250 (0.071) &&-0.175 (0.035) &&0.240 (0.053) &&-5.426 (0.057) \\
&&& [1.352; 1.644] && [0.123; 0.393] && [-0.248; -0.100] && [0.135; 0.353] && [-5.535; -5.344] \\
\cmidrule{2-2} \cmidrule{4-4} \cmidrule{6-6} \cmidrule{8-8}\cmidrule{10-10}\cmidrule{12-12}
 & \multirow{2}{*}{1.9} && 2.072 (0.109) &&0.245 (0.074) &&-0.210 (0.032) &&0.223 (0.059) &&-5.400 (0.062) \\
&&& [1.885; 2.303] && [0.068; 0.400] && [-0.273; -0.152] && [0.116; 0.349] && [-5.537; -5.244] \\
\cmidrule{2-2} \cmidrule{4-4} \cmidrule{6-6} \cmidrule{8-8}\cmidrule{10-10}\cmidrule{12-12}
 & \multirow{2}{*}{2.5} && 2.720 (0.148) &&0.308 (0.060) &&-0.200 (0.026) &&0.198 (0.042) &&-5.367 (0.058) \\
&&& [2.412; 3.016] && [0.192; 0.426] && [-0.257; -0.155] && [0.119; 0.281] && [-5.462; -5.274] \\
\cmidrule{2-2} \cmidrule{4-4} \cmidrule{6-6} \cmidrule{8-8}\cmidrule{10-10}\cmidrule{12-12}
 & \multirow{2}{*}{5.0} && 5.311 (0.356) &&0.226 (0.045) &&-0.176 (0.019) &&0.293 (0.036) &&-5.291 (0.035) \\
&&& [4.676; 6.070] && [0.137; 0.316] && [-0.215; -0.140] && [0.225; 0.368] &&[-5.346; -5.246] \\
\bottomrule
\end{tabular}}
\end{table}

From Table \ref{table:example} it is clear that, for any $\eta_i$,   with $i\in \{1,\cdots, 5\}$,  the use of  the entire or  the  thinned (thinning parameter $\mathfrak{t} = 200$)  does not yield significant improvement  in terms of parameter estimation. Not even the differences in the sample standard deviations or  in the credibility intervals, which are the statistics affected by the sample correlation, justify the  computational effort to obtain a sample of size 200,801.  The same conclusions are obtained when  considering $d_0 \in \{0.10, 0.35, 0.45\}$.  This concludes the example.
\end{example}

Different prior distributions are tested  as explained in the sequel.   Since the conditional probability density function of $I_0$ given $\boldsymbol \eta$ is difficult to obtain, in all scenarios,   it is assumed that  $g(Z_s) = 0$, for all $s \leq 1$, and it is fixed $p_{\mbox{\tiny\ensuremath{I_0}}}(\cdot|\boldsymbol  \eta) = 1$.  Moreover, since  \eqref{fieprocess} is well defined regardless the  relation among the parameters of the model, it is  assumed that
    \[
    p_{\mbox{\tiny\ensuremath{(-i)}}}(\boldsymbol \eta_{(-i)}|\eta_i) \propto \prod_{j\neq i} \pi_j(\eta_j), \quad \mbox{for any }  i\in\{1,\cdots, 5\}.
    \]

 In a first moment  the prior distributions for $\nu, d, \theta, \gamma$ and $\omega$ are selected by considering  only the basic set of information  usually available in practice. The information on each parameter and the corresponding prior selected are given in Table \ref{table:priors}.   This scenario shall be referred to as Case 1.   Table \ref{limits} presents the  mean, standard deviation,  lower and  upper limits  for the transition kernel considered  at iteration $m$ of the Gibbs sampler with Metropolis steps,  when the  prior for  $\eta_i$  in $\boldsymbol {\eta}  =  (\nu, d, \theta, \gamma, \omega)' :=  (\eta_1, \cdots, \eta_5)' $, for each  $i\in \{1,\cdots,5\}$, is defined according to Case 1.

\begin{table}[!ht]
\renewcommand{\arraystretch}{1.2}
\centering
\caption{Information available in practice for the parameter $\eta_i$  in $\boldsymbol {\eta}  =  (\nu, d, \theta, \gamma, \omega)' :=  (\eta_1, \cdots, \eta_5)' $  and the corresponding prior considered, for each  $i\in \{1,\cdots, 5\}$.}\label{table:priors}
\begin{tabular}{p{0.75\linewidth}cp{0.17\linewidth}}
\toprule
\multicolumn{1}{c}{\bf Information Available} && \multicolumn{1}{c}{\bf Prior}\\
\cmidrule{1-1}\cmidrule{3-3}

The generalized error distribution is well defined for any $\nu >0$. && \multicolumn{1}{c}{$\nu \sim \mathds{I}_{(0,\infty)}(\nu)$ *}\\
\cmidrule{1-1}\cmidrule{3-3}

Long-memory  in volatility is observed if and only if $d \in (0, 0.5)$.
This  characteristic can be detected, for example,
 through the periodogram function of the time series $\{\ln(X_t^2)\}_{t=1}^n$ \citep[see][]{LP2013}.  && \multicolumn{1}{c}{$d \sim \mathcal{U}(0,0.5)$} \\
\cmidrule{1-1}\cmidrule{3-3}

Empirical evidence suggests that $\theta \in [-1, 0]$. ** && \multicolumn{1}{c}{$\theta \sim \mathcal{U}(-1,0)$}  \\

\cmidrule{1-1}\cmidrule{3-3}

 Empirical evidence suggests that $\gamma \in [0, 1]$. ** && \multicolumn{1}{c}{$\gamma \sim \mathcal{U}(0,1)$} \\
\cmidrule{1-1}\cmidrule{3-3}

 $\omega = \mathds{E}(\ln(h_t^2)) = \mathds{E}(\ln(X_t^2)) + \mathds{E}(\ln(Z_t^2))$.\newline
The choice of the interval for $\omega$ will depend on the magnitude of the data.
The sample mean of $\{\ln(X_t^2)\}_{t=1}^n$ or $\ln(\hat \sigma_X^2)$, where
  $\hat \sigma_X^2$ is the sample variance of $\{X_t\}_{t=1}^n$, can be used to obtain a rough approximation for $\omega$ && \multicolumn{1}{c}{$\omega \sim \mathcal{U}(-15,15)$.}\\
\bottomrule
\end{tabular}
\begin{description}
\item[{\footnotesize {\bf Notes:}}]{\footnotesize *  Given  $A \subset \mathds{R}$, the symbol $\mathds{I}_{A}(x)$  denotes the improper prior defined as   $1$,  if  $x\in A$, and 0, if  $x\notin A$.\\
 **  See, for instance, \cite{NE1991, BM1996, RV2008, LP2013}.  To the best of our  knowledge,   a FIEGARCH model for which  $\theta$ or $\gamma$ are not in the intervals, respectively,  $[-1,0]$ and  $[0, 1]$  has never been reported in the literature. }
\end{description}
\end{table}

\begin{table}[!ht]
\centering
\caption{Parameters of the truncated normal distribution (transition kernel) considered,   at iteration $m$ of the Gibbs sampler, to obtain the sample from the posterior distribution of  the parameter $\eta_i$  in $\boldsymbol {\eta} =  (\nu, d, \theta, \gamma, \omega)' :=  (\eta_1, \cdots, \eta_5)' $, for each  $i\in \{1,\cdots, 5\}$.}\label{limits}
\begin{tabular*}{1\textwidth}{@{\extracolsep{\fill}} crrrrr}
\toprule
Parameter & $\nu$ & $d$ & $\theta$ & $\gamma$ & $\omega$\\
\midrule
Mean ($y$)  & $\nu^{(m-1)}$ & $d^{(m-1)}$ & $\theta^{(m-1)}$ & $\gamma^{(m-1)}$ & $\omega^{(m-1)}$\\
Standard Deviation ($\sigma$) &  $0.500$  & $0.025$ & $0.050$ &   $0.050$  & $1.500$ \\
Lower Limit ($a$) & $\phantom{-}0.000$   & $0.000$ & $-1.000$ & $0.000$ & $-15.000$ \\
Upper Limit ($b$) & $10.000$    & $0.500$ & $0.000$ & $1.000$& $15.000$ \\
\bottomrule
\end{tabular*}\vspace{-0.3cm}\\
\begin{description}
\item[{\footnotesize {\bf Note:}}]{\footnotesize $\eta_i^{(m-1)}$, for any $i\in \{1,\cdots, 5\}$,  denotes the parameter value obtained in the $(m-1)$th  iteration.   Different combinations of  standard deviation,  lower and upper limits  were tested for the parameter $\eta_i$ in $\boldsymbol {\eta} =  (\nu, d, \theta, \gamma, \omega)'$, for each $i\in\{1,\cdots, 5\}$. The values presented here correspond to the final choice.}
\end{description}
\end{table}

 In a second moment  the knowledge  on the true parameter values is gradually incorporated to provide  more informative priors for $d, \theta$ and/or $\gamma$.   This analysis, combined with the first scenario, provides information on the sensitivity of the estimates with respect to the priors functions and hyperparameters.  In all cases, the  priors for $\nu$ and $\omega$ are the same and are the ones defined in Table \ref{table:priors}.   The scenarios considered in this second step are described in the following and shall be referred to as Case 2 -  Case 5.

\subsubsection*{Case 2:  Gaussian Prior for $x = \phi^{-1}(d)$ and Uniform Priors for $\theta$ and $\gamma$.}

In this case $\theta$ and $\gamma$ remain with the same priors as in Case 1. For the parameter  $d$  it is  assumed that    $x \sim \mathcal{N}(\mu_\phi,\sigma_\phi^2)$ and   $d = \phi(x)$, where $\phi: \mathbb{R} \to (0,0.5)$  is given by
    \begin{equation}\label{eq:phi}
    \phi(x) = \frac{e^x}{2(1+ e^x)}, \quad \mbox{for all }  x\in\mathds{R}.
    \end{equation}

First, the knowledge of  $d_0$ is applied to set $\mu_{\phi} = \phi^{-1}(d_0)$, so  $\mu_{\phi} \in \{ -1.386,  0.000,  0.847,  2.197\}$, respectively, for $d_0 \in\{0.10, 0.25, 0.35, 0.45\}$.  This scenario shall be referred to as  C2.1.  Second,  the  knowledge on $d_0$ is ignored and the parameter $\mu_\phi$ is assumed to be equal to zero.  This scenario shall be referred to as  C2.2.  For both,   C2.1 and  C2.2,   different values of $\sigma_\phi$ are tested.  Third,  the approaches considered in   C2.1 and  C2.2  are combined  by setting  $\mu_{\phi} = \phi^{-1}(\bar d)$, where $\bar d$ is the estimate of $d$ obtained in  C2.2.    This scenario shall be referred to as  C2.3.  The value of  $\sigma_\phi$ considered in   C2.3 is the one  which  provides better estimates for $d$  in  C2.1.

 The kernel parameter values for $\nu, \theta, \gamma$ and $\omega$ are the same  as in the  Case 1.   For $x = \phi^{-1}(d)$,  at
 iteration $m$ of the Gibbs sampler, the kernels mean ($y$),   standard deviation ($sd$), lower ($a$) and upper limits ($b$) are set, respectively, as  $x^{(m-1)}$,  1, -10 and 10, where  $x^{(m-1)}$ is the parameter value obtained at iteration $m-1$.

\subsubsection*{Case 3:   Gaussian Prior for $x = \phi^{-1}(d)$,  Beta Prior for $-\theta$ and Uniform Prior for $\gamma$.}

In this case,  the priors of   $\gamma$ and $d$ are  the same ones considered, respectively,   in  Case 1 and  in  scenario  C2.1 of  Case 2.  It is  also assumed that  $-\theta \sim \mbox{Beta}(a_1,b_1)$, which is equivalent to set
\[
\pi_3(\theta) = (-\theta)^{a_1-1}(1+\theta)^{b_1-1}B(a_1,b_1), \quad \theta \in [-1, 0],
\]
where $B(\cdot, \cdot)$ is the beta function.

 First, the fact that   $X \sim \mbox{Beta}(a,b)$  implies $\mathbb{E}(X) = a(a+b)^{-1}$, is  applied to set  $b_1 = a_1(1+\theta_0)(-\theta_0)^{-1}$, where $\theta_0 = -0.15$ is  the true parameter value considered in this simulation study.   Different values of $a_1$ are tested.   This scenario shall be referred to as  C3.1.  Second,   the  knowledge on $\theta_0$ is ignored and different combinations of $a_1$ and $b_1$ are tested.  This scenario shall be referred to as  C3.2.   Third,  the approaches considered in   C3.1 and  C3.2  are combined  by setting  $b_1 = a_1(1+\bar\theta_0)(-\bar\theta_0)^{-1}$,  where $-\bar \theta_0$ is the estimate of $\theta$ obtained in  C3.2.    The value of  $a_1$ considered in  this case is the one  which  provides better estimates for $\theta$  in  C3.1.   This scenario shall be referred to as  C3.3.

 The kernel parameter values are the same  as  in   Case 2.

\subsubsection*{Case 4:  Gaussian Prior for $x = \phi^{-1}(d)$, Beta Priors for $-\theta$ and $\gamma$.}

In this case,  the priors of   $d$ and $-\theta$ are  the same ones considered, respectively,    in  scenario  C2.1 of  Case 2 and in scenario  C3.1 of Case 3.   It is also  assumed that  $\gamma \sim \mbox{Beta}(a_2,b_2)$.  Two  scenarios, denoted by  C4.1 and   C4.2 are considered.  With the obvious identifications, the construction of  C4.1 and  C4.2 is analogous, respectively,  to the construction of scenarios  C3.1 and  C3.2 in  Case 3.

 The kernel parameter values are the same  as  in   Case 2.

\subsubsection*{Case 5:   Beta Priors for $2d$, $-\theta$ and $\gamma$.}

In this case,  the priors of  $-\theta$ and $\gamma$ are  the same ones considered,  respectively,    in  scenario  C3.1 of  Case 3 and in scenario  C4.1 of  Case 4.  Moreover, for each  $d_0 \in \{0.10, 0.25, 0.35, 0.45\}$ considered in this simulation study,   it is assumed that
 $2d \sim  \mbox{Beta}(a_3,b_3)$, which is equivalent to set
 \[
\pi_2(d) =  2(2d)^{a_3-1}(1- 2d)^{b_3-1}B(a_3,b_3), \quad d \in [0, 0.5],
\]
where $B(\cdot, \cdot)$ is the beta function.

In this case, only two scenarios are considered. First, it is assumed that  $b_3 = a_3(1-2d_0)(2d_0)^{-1}$ and different values of $a_3$ are tested. This scenario shall be referred to as  C5.1.   Second,   an approach similar to scenarios  C3.3 and  C4.3, respectively,  in  Case 3 and   Case 4, is considered.   However, in this case,  it is assumed that $b_3 = a_3(1-2\bar d)(2\bar d)^{-1}$, with $\bar d$ obtained in  Case 1.  The value of  $a_3$ considered in  this case is the one  which  provides better estimates for $d$  in  C5.1.   This scenario shall be referred to as  C5.2.

The kernel parameter values are the same  as  in   Case 1.

\subsection{Estimates and Performance Measures}

Let $\{\eta_i^{(k)}\}_{k=1}^M$ be a sample of size $M$ from the posteriori distribution of  $\eta_i$  in $\boldsymbol {\eta} =  (\nu, d, \theta, \gamma, \omega)' :=  (\eta_1, \cdots, \eta_5)' $, for any  $i\in \{1,\cdots, 5\}$. Denote by $\bar \eta_i$ and $\mbox{sd}_{\eta_i}$,  respectively,  the sample mean and standard deviation of $\{\eta_i^{(k)}\}_{k=1}^M$, namely,
    \[
     \bar \eta_i = \frac{1}{M}\sum_{k=1}^M\eta_i^{(k)} \quad \mbox{and} \quad \mbox{sd}_{\eta_i} =\sqrt{\frac{1}{M}\sum_{k=1}^M(\eta_i^{(k)} - \bar\eta_i)^2}, \quad \mbox{for any} \quad i\in \{1,\cdots, 5\}.
    \]
Then  the  estimate  $\hat \eta_i$ of $\eta_i$ is defined as $\hat \eta_i := \bar \eta_i$.

Moreover,  let  $\hat q_i(\alpha)$  denote the  $\alpha$ quantile\footnote{In this work, the following definition is adopted \citep{BD1991}.  Given any $ 0 \leq  \alpha \leq 1$, the number $q(\alpha)$ satisfying $\mathds{P}(X \leq q(\alpha)) \geq \alpha$ and $\mathds{P}(X \geq q(\alpha)) \geq 1- \alpha$,
is called a quantile of order $\alpha$ (or $\alpha$ quantile) for the random variable $X$ (or for the distribution function of $X$).  }     for the  posterior sample distribution of $\eta_i$, for any $\alpha \in [0,1]$ and $i\in \{1,\cdots, 5\}$.   Then   a $100(1-\alpha)\%$  credibility interval for $\eta_i$ is given by
\[
CI_{1-\alpha}(\eta_i) = \Big[ \hat q_i\Big(\frac{\alpha}{2}\Big), \hat q_i\Big(\frac{1-\alpha}{2}\Big) \Big], \quad \mbox{for any}\quad  i\in \{1,\cdots, 5\}.
\]

Furthermore, the estimation bias  and the  absolute percentage error (ape) of estimation  are given, respectively,   by
    \[
    \mbox{bias}_{\eta_i} = \bar \eta_i - \eta_i \quad  \mbox{and} \quad  \mbox{ape}_{\eta_i} = \bigg|\frac{\mbox{bias}_{\eta_i}}{\eta_i}\bigg|,  \quad \mbox{for any} \quad i\in \{1,\cdots, 5\}.
    \]

\subsection{Results}

The results obtained in this simulation study, by considering the scenarios described in Section \ref{section:priors}, are the following.

\subsubsection*{Case 1:  The Priors as Defined in Table \ref{table:priors}.}

Table \ref{table:uniform} present the  summary statistics  for the samples obtained from  the posterior distribution  for  each parameter of the  FIEGARCH$(0,d,0)$ model. The statistics reported in this table (the same applies to Table \ref{table:beta}) are the sample  mean ($\bar \eta_i$), the sample standard deviation ($\mbox{sd}_{\eta_i}$) and the 95\% credibility interval $CI_{0.95}(\eta_i)$ for the parameter $\eta_i$  in $\boldsymbol {\eta} =  (\nu, d, \theta, \gamma, \omega)' :=  (\eta_1, \cdots, \eta_5)' $, for each  $i\in \{1,\cdots, 5\}$.    The bold-face font for the mean indicates that the absolute percentage error of estimation ($\mbox{ape}_{\eta_i}$) in the corresponding case is higher than 0.10 (that is, 10\%). The bold-face font for the credibility interval indicates that the true parameter value is not contained in the interval.

\begin{table}[!htp]
\centering
\setlength{\tabcolsep}{5pt}
\caption{Summary for the sample obtained from posterior distributions considering all prior uniforms:  mean $\bar \eta_i$, standard deviation $\mbox{sd}_{\eta_i}$ and the 95\% credibility interval $CI_{0.95}(\eta_i)$ for the parameter $\eta_i$  in $\boldsymbol {\eta} =  (\nu, d, \theta, \gamma, \omega)' :=  (\eta_1, \cdots, \eta_5)' $, for each  $i\in \{1,\cdots, 5\}$.  The true parameter values considered in this simulation are $d_0 \in\{0.10, 0.25, 0.35, 0.45\}$, $\nu_0 \in \{1.1, 1.5, 2.5, 5.0\}$, $\theta_0 = -0.15$, $\gamma_0 = -0.24$ and $\omega_0 = -5.4$.}\label{table:uniform}
{\footnotesize
\begin{tabular}{c|ccccccccccc}
\toprule
\multicolumn{1}{c}{\multirow{2}{*}{$d_0$}}&
\multicolumn{1}{c}{\multirow{2}{*}{$\nu_0$}}&&
\multicolumn{1}{c}{$\bar{\nu}$ $(\mbox{sd}_{\nu})$} &&
\multicolumn{1}{c}{$\bar{d}$ $(\mbox{sd}_{d})$} &&
\multicolumn{1}{c}{$\bar{\theta}$ $(\mbox{sd}_{\theta})$} &&
\multicolumn{1}{c}{$\bar{\gamma}$ $(\mbox{sd}_{\gamma})$} &&
\multicolumn{1}{c}{$\bar{\omega}$ $(\mbox{sd}_{\omega})$} \\
\multicolumn{1}{c}{}&\multicolumn{1}{c}{}&&
\multicolumn{1}{c}{$CI_{0.95}(\nu)$} &&
\multicolumn{1}{c}{$CI_{0.95}(d)$} &&
\multicolumn{1}{c}{$CI_{0.95}(\theta)$} &&
\multicolumn{1}{c}{$CI_{0.95}(\gamma)$} &&
\multicolumn{1}{c}{$CI_{0.95}(\omega)$} \\

\cmidrule{1-2} \cmidrule{4-4} \cmidrule{6-6} \cmidrule{8-8}\cmidrule{10-10}\cmidrule{12-12}
\multirow{10}{*}{0.10}
 & \multirow{2}{*}{1.1} && 1.093 (0.044) &&\textbf{0.181} (0.123) &&\textbf{-0.084} (0.041) &&0.236 (0.066) &&-5.438 (0.058) \\
&&& [0.989; 1.195] && [0.005; 0.458] && [-0.171; -0.013] && [0.093; 0.357] && [-5.551; -5.337] \\
\cmidrule{2-2} \cmidrule{4-4} \cmidrule{6-6} \cmidrule{8-8}\cmidrule{10-10}\cmidrule{12-12}
 & \multirow{2}{*}{1.5} && 1.480 (0.069) &&\textbf{0.147} (0.079) &&\textbf{-0.177} (0.038) &&0.232 (0.052) &&-5.420 (0.036) \\
&&& [1.353; 1.635] && [0.020; 0.330] && [-0.258; -0.106] && [0.122; 0.340] && [-5.510; -5.338] \\
\cmidrule{2-2} \cmidrule{4-4} \cmidrule{6-6} \cmidrule{8-8}\cmidrule{10-10}\cmidrule{12-12}
 & \multirow{2}{*}{1.9} && 2.088 (0.111) &&0.093 (0.055) &&\textbf{-0.220} (0.032) &&\textbf{0.216} (0.060) &&-5.410 (0.035) \\
&&& \textbf{[1.908; 2.296]} && [0.004; 0.201] && \textbf{[-0.286; -0.154]} && [0.105; 0.337] && [-5.486; -5.340] \\
\cmidrule{2-2} \cmidrule{4-4} \cmidrule{6-6} \cmidrule{8-8}\cmidrule{10-10}\cmidrule{12-12}
 & \multirow{2}{*}{2.5} && 2.724 (0.140) &&\textbf{0.192} (0.076) &&\textbf{-0.201} (0.025) &&\textbf{0.198} (0.045) &&-5.388 (0.031) \\
&&& [2.491; 3.027] && [0.038; 0.330] && \textbf{[-0.256; -0.153]} && [0.116; 0.287] && [-5.448; -5.333] \\
\cmidrule{2-2} \cmidrule{4-4} \cmidrule{6-6} \cmidrule{8-8}\cmidrule{10-10}\cmidrule{12-12}
 & \multirow{2}{*}{5.0} && 5.297 (0.364) &&0.101 (0.051) &&\textbf{-0.174} (0.020) &&\textbf{0.297} (0.036) &&-5.336 (0.028) \\
&&& [4.641; 6.014] && [0.015; 0.217] && [-0.215; -0.133] && [0.232; 0.374] && \textbf{[-5.383; -5.287}] \\
\midrule

\multirow{10}{*}{0.25}
 & \multirow{2}{*}{1.1} && 1.108 (0.039) &&0.265 (0.129) &&\textbf{-0.089} (0.038) &&0.230 (0.060) &&-5.476 (0.052) \\
&&& [1.028; 1.198] && [0.016; 0.467] && [-0.176; -0.019] && [0.108; 0.345] && [-5.553; -5.366] \\
\cmidrule{2-2} \cmidrule{4-4} \cmidrule{6-6} \cmidrule{8-8}\cmidrule{10-10}\cmidrule{12-12}
 & \multirow{2}{*}{1.5} && 1.474 (0.067) &&0.250 (0.071) &&\textbf{-0.175} (0.035) &&0.240 (0.053) &&-5.426 (0.057) \\
&&& [1.352; 1.644] && [0.123; 0.393] && [-0.248; -0.100] && [0.135; 0.353] && [-5.535; -5.344] \\
\cmidrule{2-2} \cmidrule{4-4} \cmidrule{6-6} \cmidrule{8-8}\cmidrule{10-10}\cmidrule{12-12}
 & \multirow{2}{*}{1.9} && 2.072 (0.109) &&0.245 (0.074) &&\textbf{-0.210} (0.032) &&0.223 (0.059) &&-5.400 (0.062) \\
&&& [1.885; 2.303] && [0.068; 0.400] && \textbf{[-0.273; -0.152]} && [0.116; 0.349] && [-5.537; -5.244] \\
\cmidrule{2-2} \cmidrule{4-4} \cmidrule{6-6} \cmidrule{8-8}\cmidrule{10-10}\cmidrule{12-12}
 & \multirow{2}{*}{2.5} && 2.720 (0.148) &&\textbf{0.308} (0.060) &&\textbf{-0.200} (0.026) &&\textbf{0.198} (0.042) &&-5.367 (0.058) \\
&&& [2.412; 3.016] && [0.192; 0.426] && \textbf{[-0.257; -0.155]} && [0.119; 0.281] && [-5.462; -5.274] \\
\cmidrule{2-2} \cmidrule{4-4} \cmidrule{6-6} \cmidrule{8-8}\cmidrule{10-10}\cmidrule{12-12}
 & \multirow{2}{*}{5.0} && 5.311 (0.356) &&0.226 (0.045) &&\textbf{-0.176} (0.019) &&\textbf{0.293} (0.036) &&-5.291 (0.035) \\
&&& [4.676; 6.070] && [0.137; 0.316] && [-0.215; -0.140] && [0.225; 0.368] && \textbf{[-5.346; -5.246]} \\
\midrule

\multirow{10}{*}{0.35}
 & \multirow{2}{*}{1.1} && 1.099 (0.040) &&0.349 (0.108) &&\textbf{-0.097} (0.038) &&0.230 (0.056) &&-5.495 (0.090) \\
&&& [1.009; 1.194] && [0.093; 0.492] && [-0.178; -0.027] && [0.121; 0.330] && [-5.674; -5.318] \\
\cmidrule{2-2} \cmidrule{4-4} \cmidrule{6-6} \cmidrule{8-8}\cmidrule{10-10}\cmidrule{12-12}
 & \multirow{2}{*}{1.5} && 1.479 (0.065) &&0.329 (0.065) &&\textbf{-0.178} (0.036) &&0.246 (0.052) &&-5.423 (0.076) \\
&&& [1.352; 1.639] && [0.204; 0.461] && [-0.246; -0.106] && [0.143; 0.340] && [-5.561; -5.302] \\
\cmidrule{2-2} \cmidrule{4-4} \cmidrule{6-6} \cmidrule{8-8}\cmidrule{10-10}\cmidrule{12-12}
 & \multirow{2}{*}{1.9} && 2.064 (0.110) &&0.364 (0.054) &&\textbf{-0.199} (0.030) &&0.233 (0.050) &&-5.377 (0.090) \\
&&& [1.843; 2.299] && [0.227; 0.461] && [-0.265; -0.139] && [0.139; 0.330] && [-5.535; -5.199] \\
\cmidrule{2-2} \cmidrule{4-4} \cmidrule{6-6} \cmidrule{8-8}\cmidrule{10-10}\cmidrule{12-12}
 & \multirow{2}{*}{2.5} && 2.732 (0.150) &&0.380 (0.052) &&\textbf{-0.201} (0.024) &&\textbf{0.200} (0.043) &&-5.307 (0.066) \\
&&& [2.481; 3.031] && [0.283; 0.479] && \textbf{[-0.254; -0.161]} && [0.110; 0.285] && [-5.410; -5.149] \\
\cmidrule{2-2} \cmidrule{4-4} \cmidrule{6-6} \cmidrule{8-8}\cmidrule{10-10}\cmidrule{12-12}
 & \multirow{2}{*}{5.0} && 5.229 (0.321) &&0.318 (0.040) &&\textbf{-0.177} (0.019) &&\textbf{0.289} (0.036) &&-5.227 (0.046) \\
&&& [4.603; 5.864] && [0.243; 0.409] && [-0.216; -0.140] && [0.227; 0.366] && \textbf{[-5.298; -5.127]} \\
\midrule

\multirow{10}{*}{0.45}
 & \multirow{2}{*}{1.1} && 1.096 (0.039) &&0.436 (0.053) &&\textbf{-0.115} (0.034) &&0.241 (0.054) &&-5.453 (0.128) \\
&&& [1.010; 1.161] && [0.313; 0.499] && [-0.187; -0.047] && [0.136; 0.338] && [-5.716; -5.174] \\
\cmidrule{2-2} \cmidrule{4-4} \cmidrule{6-6} \cmidrule{8-8}\cmidrule{10-10}\cmidrule{12-12}
 & \multirow{2}{*}{1.5} && 1.475 (0.073) &&0.411 (0.048) &&\textbf{-0.179} (0.034) &&0.257 (0.048) &&-5.411 (0.123) \\
&&& [1.353; 1.627] && [0.311; 0.494] && [-0.246; -0.110] && [0.158; 0.353] && [-5.600; -5.133] \\
\cmidrule{2-2} \cmidrule{4-4} \cmidrule{6-6} \cmidrule{8-8}\cmidrule{10-10}\cmidrule{12-12}
 & \multirow{2}{*}{1.9} && 2.052 (0.111) &&0.450 (0.032) &&\textbf{-0.191} (0.026) &&0.243 (0.043) &&-5.367 (0.130) \\
&&& [1.846; 2.279] && [0.385; 0.497] && [-0.238; -0.141] && [0.165; 0.320] && [-5.614; -5.125] \\
\cmidrule{2-2} \cmidrule{4-4} \cmidrule{6-6} \cmidrule{8-8}\cmidrule{10-10}\cmidrule{12-12}
 & \multirow{2}{*}{2.5} && 2.725 (0.152) &&0.447 (0.032) &&\textbf{-0.206} (0.021) &&\textbf{0.211} (0.041) &&-5.150 (0.081) \\
&&& [2.461; 3.021] && [0.381; 0.495] && \textbf{[-0.247; -0.167]} && [0.133; 0.296] && \textbf{[-5.310; -4.985]} \\
\cmidrule{2-2} \cmidrule{4-4} \cmidrule{6-6} \cmidrule{8-8}\cmidrule{10-10}\cmidrule{12-12}
 & \multirow{2}{*}{5.0} && 5.177 (0.322) &&0.417 (0.032) &&\textbf{-0.177} (0.019) &&\textbf{0.286} (0.032) &&-5.041 (0.068) \\
&&& [4.553; 5.832] && [0.348; 0.480] && [-0.220; -0.140] && [0.228; 0.350] && \textbf{[-5.182; -4.929]} \\
\bottomrule
\end{tabular}}
\begin{description}
\item[{\footnotesize {\bf Note:}}]{\footnotesize  The  bold-face font for the estimated \textbf{mean} indicates that the absolute percentage error is higher than 10\%.  The bold-face font for the \textbf{credibility interval} indicates that the interval does not contain the true parameter value.}
\end{description}
\end{table}

From  Table \ref{table:uniform} one observes  that  the parameters $\nu$ and $\omega$  are always well estimated, in terms of absolute percentage error (ape),  regardless the combination of $d_0 \in\{0.10, 0.25, 0.35, 0.45\}$ and   $\nu_0 \in \{1.1, 1.5, 2.5, 5.0\}$ considered (the error is less than 10\% in all cases).     The credibility interval $CI_{0.95}(\nu)$ contains the true parameter value $\nu_0$ in all cases, except when $d_0 = 0.10$ and $\nu_0 = 1.9$.  Also, the estimation bias for $\nu$ is always negative when $\nu_0 < 1.9$  and positive when $\nu_0 \geq 1.9$, except for the combination $(\nu_0, d_0) = (1.1, 0.25)$.    For the parameter $\omega$,  the credibility interval $CI_{0.95}(\omega)$ does not contain the true parameter value ($\omega_0 = -5.4$) in 5 out of 20 combinations of $\nu_0$ and $d_0$ (see $\nu_0 = 5$ and all $d_0$; $\nu_0 = 2.5$ and $d_0 = 0.45$).  Moreover,   the estimation  bias for $\omega$ is always negative when $\nu_0 \leq 1.5$ and always positive when $\nu_0 \geq 2.5$.

 Table  \ref{table:uniform}  also reports that $\mbox{ape}_\theta >  10\%$ for all combinations of $d_0$ and $\nu_0$.  On the other hand, in most cases (14 out of 20),   the credibility interval $CI_{0.95}(\theta)$ contains the true parameter value $\theta_0 = -0.15$.   The cases for which  $\theta_0 \notin CI_{0.95}(\theta)$ are   $\nu = 1.9$ and $d_0 \in \{0.10, 0.25\}$ and  $\nu_0 = 2.5$ and any $d_0$.   The bias for $\theta$ is always positive when $\nu_0 = 1.1$ (for any $d_0$) and negative in all other cases.

 Furthermore, Table \ref{table:uniform} shows that the parameter $\gamma$ seems to be better estimated when the GED distribution presents heavy tails ($\nu_0 < 2$), except when $d = 0.10$, in which case $\mbox{ape}_\gamma >  10\%$  when $\nu_0 = 1.9$.  Also, with exception of four cases ($d_0 = 0.10$ and $\nu_0 \in \{1.1, 1.5, 2.5\}$; $d_0  = 0.25$ and $\nu_0 = 2.5$), the parameter $d$ is always well estimated.  The bias for  parameters $d$ and $\gamma$ does not seem to follow any pattern and both, $d_0 \in CI_{0.95}(d)$ and $\gamma_0 \in CI_{0.95}(\gamma)$, for any combination of $\nu_0$ and $d_0$.

\subsubsection*{Case 2:  Gaussian Prior for $x = \phi^{-1}(d)$ and Uniform Priors for $\theta$ and $\gamma$.}

Changing the prior for $d$ does not yield  significant difference on the  estimation of $\nu$, $\theta$, $\gamma$ and  $\omega$.

When the true  value of  $d_0$ is used to set $\mu_\phi = \phi^{-1}(d_0)$ (scenario  C2.1),  the  best performance is observed by letting $\sigma_\phi = 0.15$.  In this case,  the absolute percentage error of estimation ($\mbox{ape}_d$) is smaller than 10\% for all combinations of $\nu_0$ and $d_0$,  with $d_0 \in\{0.10, 0.25, 0.35, 0.45\}$ and   $\nu_0 \in \{1.1, 1.5, 2.5, 5.0\}$.   If $\sigma_{\phi} = 0.10$ the chain takes too long to move from the initial point when $d_0 = 0.10$.   When $\sigma_\phi = 0.25$,  there is  only one case for which $\mbox{ape}_d  > 10\%$  ($d_0 = 0.10$ and $\nu_0 = 2.5$).  In fact, in this case, $\mbox{ape}_d = 0.103$, which is still acceptable ($\sigma_\phi = 0.15$  still seems to be the best choice).  Furthermore,  as $\sigma_\phi$ increases, the number of cases for which $\mbox{ape}_d  > 10\%$ also increases. For instance, when $\sigma_\phi \in \{0.50, 1.00, 3.00\}$,  $\mbox{ape}_d  >  10\%$ in  2, 4 and 10 cases, respectively.

When $d_0$ is assumed unknown and $\mu_\phi$ is set to zero (scenario  C2.2),   $\sigma_\phi = 3$ seems to  provide better results than smaller values of $\sigma_\phi$.    Under this scenario,   $\mbox{ape}_d  > 10\%$  for  8 out of 20 combinations of $\nu_0$ and $d_0$. Therefore,   $d\sim \mathcal{U}(0,0.5)$ still provides better estimates for the parameter $d$ (see Table \ref{table:uniform}).   Higher values of $\sigma_\phi $ do not improve the  estimation of $d$.  Too high values of $\sigma_\phi$ actually  make  the estimation worst.   In particular, when $\sigma_\phi  = 4$   the results are similar to $\sigma_\phi  = 3$, if $d_0 > 0.1$.   If $d_0 = 0.1$ then $\sigma_\phi  = 3$ is slightly better than $\sigma_\phi  = 4$.  When $\sigma_\phi  = 5$,   $\mbox{ape}_d$ is, in most cases,  higher than when $\sigma_\phi  = 3$.     For $\sigma_\phi $ smaller than 3 the estimation bias is much higher.    For instance,  when $\sigma_\phi = 0.15$,   $\mbox{ape}_d \leq 10\%$ only for $d_0 = 0.25$ (for all $\nu_0$).  For all other combinations of $d_0$ and $\nu_0$ $\mbox{ape}_d >  20\%$.    Also, when $b = 1$,      $\mbox{ape}_d > 20\%$   in 12 out of 20 cases. In particular,   $\mbox{ape}_d > 20\%$  for $d_0 = 0.10$ and  all $\nu_0$.   As it should be expected,  C2.1 performs much better than  C2.2.

Upon considering a two step estimator (scenario  C2.3),  no improvement is observed, when compared to scenario  C2.2. In fact, the  estimates obtained by letting $\mu_{\phi} = \phi^{-1}(\bar d)$ (where $\bar d$ is the estimate of $d$ obtained in  C2.2) and  $\sigma_\phi = 0.15$ (the parameter which leads to the best performance in  C2.1)  are very close to $\bar d$ itself.

\subsubsection*{Case 3:  Gaussian Prior for $x = \phi^{-1}(d)$, Beta Prior for $-\theta$ and Uniform Prior for $\gamma$.}

The  estimation of $\nu$, $d$, $\gamma$ and  $\omega$ is not significantly affected by the change in the prior for $-\theta$.

When  the knowledge on the true parameter values $d_0$ and $\theta_0$ is applied to set  $\mu_\phi = \phi^{-1}(d_0)$,  $\sigma_\phi = 0.15$,  for each $d_0 \in\{0.10, 0.25, 0.35, 0.45\}$  (the best scenario in  Case 2), and $b_1 = a_1(1+\theta_0)(-\theta_0)^{-1}$ (scenario  C3.1), it is observed that larger values of $a_1$ lead to better estimates for $\theta$.   Although  any  $a_1 \in \{110, 150, 200\}$ leads to $\mbox{ape}_\theta \leq 10\%$, for all combinations of $d_0 \in\{0.10, 0.25, 0.35, 0.45\}$ and   $\nu_0 \in \{1.1, 1.5, 2.5, 5.0\}$,  the best performance is obtained by setting $a_1 = 110$ ($b_1 \approx 623.333$).  As $a_1$ decreases, the estimation performance decreases.  For instance,  when $a_1 = 100$, one case for which $\mbox{ape}_\theta > 10\%$ is observed. When    $a_1 = 20$  the number of cases for which $\mbox{ape}_\theta > 10\%$ increases to  10 and no case for which $\mbox{ape}_\theta \leq 10\%$ is observed if $a_1 \in \{ 2.0, 0.1\}$.   More specifically, for any $a_1 \in \{ 2.0, 0.1\}$,  $10\% <  \mbox{ape}_\theta  \leq 20\%$  for $\nu_0 \in\{1.5, 5.0\} $ and all values of $d_0$ (8 out of 20 cases) and, in the remaining  12 cases, $\mbox{ape}_\theta >  20\%$.

By assuming $\theta_0$ unknown (scenario  C3.2) or by considering a two step estimator (scenario  C3.3),  no case for which  $\mbox{ape}_\theta <  10\%$ is observed.  The combinations of $a_1$ and $b_1$ tested in scenario  C3.2  are:  $(a_1,b_1) \in \{$(2,\,3), (2,\,5), (2,\,9), (4,\,7), (5,\,7), (10,\,40), (10,\,60), (10,\,70), (100,\,500), (100,\,600)$\}$.   Among these values, the best performance is obtained when $a_1 = 10$ and $b_1 = 50$. In this case, $10\% < \mbox{ape}_\theta \leq  20\%$ in 12 out of 20 cases, which is slightly  better than the performance obtained assuming $\theta \sim \mathcal{U}(0,1)$ (in this case, $10\% < \mbox{ape}_\theta \leq  20\%$ in 8 out of 20 cases).

\subsubsection*{Case 4:  Gaussian Prior for $x = \phi^{-1}(d)$, Beta Priors for $-\theta$ and  $\gamma$.}

Analogously to  Case 2 and   Case 3, the  estimation of $\nu$, $d$, $\theta$ and  $\omega$ is not significantly affected by the change in the prior for $\gamma$.

\begin{figure}[!htb]
\includegraphics[width = 1\linewidth]{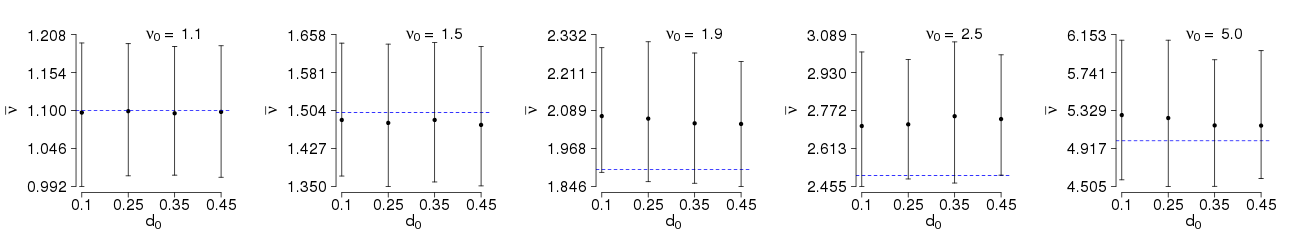}
\includegraphics[width = 1\linewidth]{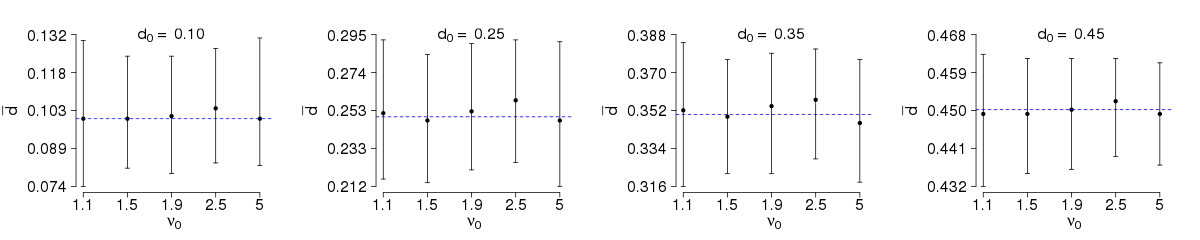}
\includegraphics[width = 1\linewidth]{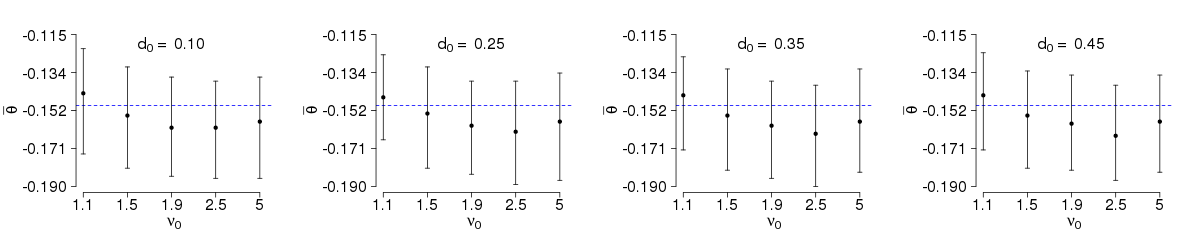}
\includegraphics[width = 1\linewidth]{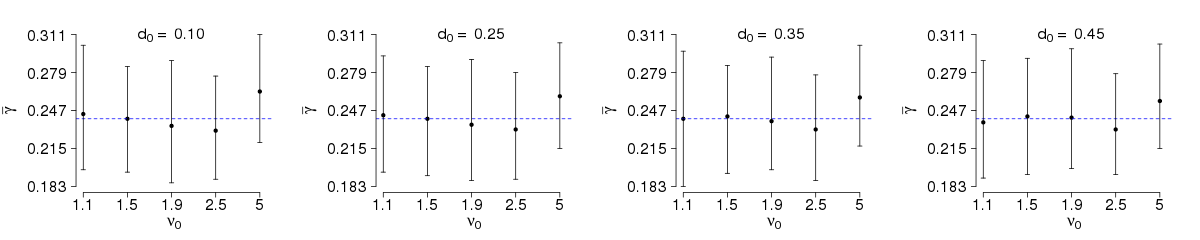}
\includegraphics[width = 1\linewidth]{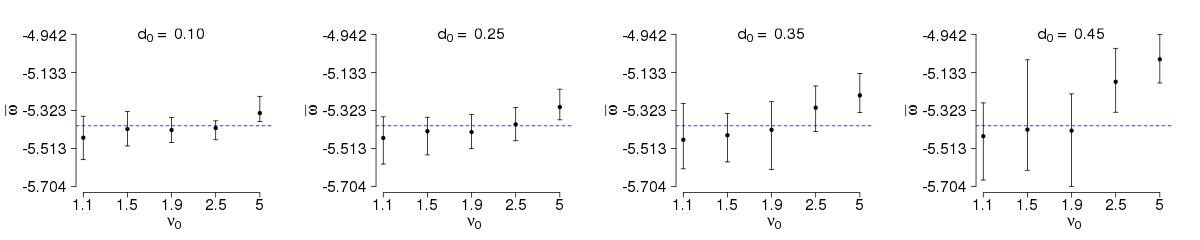}
\caption{Posterior mean  (solid circle), the true parameter value (dashed line) and the 95\% credibility interval  (solid line) for the parameters $\nu, d, \theta, \gamma$ and  $\omega$ (from top to bottom), for each combination of $d_0$ and  $\nu_0$.    The posterior distributions were obtained by  considering an  improper prior for $\nu$, a   Gaussian prior for $\phi^{-1}(d)$, Beta priors for $-\theta$ and $\gamma$ and a uniform prior for $\omega$.   The true parameters values considered in this simulation are $d_0 \in\{0.10, 0.25, 0.35, 0.45\}$, $\nu_0 \in \{1.1, 1.5, 2.5, 5.0\}$, $\theta_0 = -0.15$, $\gamma_0 = -0.24$ and $\omega_0 = -5.4$.}\label{table:gaussian}
\end{figure}

By considering the true parameter values $d_0$,  $\theta_0$ and $\gamma_0$ and setting  $\mu_\phi = \phi^{-1}(d_0)$,  $\sigma_\phi = 0.15$,  for each $d_0 \in\{0.10, 0.25, 0.35, 0.45\}$  (the best scenario in  Case 2),  $a_1 = 110$,  $b_1 = a_1(1+\theta_0)(-\theta_0)^{-1}$ (the best scenario in  Case 3) and  $b_2 = a_2(1 - \gamma_0)\gamma_0^{-1} $  (scenario  C4.1), it is observed the following:    larger values of $a_2$  (smaller than $a_1$, however) lead to better estimates for $\gamma$ and  as $a_2$ decreases, the estimation performance decays.  For instance, when $a_2 = 40$  only one case for which $\mbox{ape}_\theta > 10\%$ is observed and when   $a_2 \in \{10, 25, 30\}$,  the number of cases increases to 5 ($d_0 = 0.45$ and all $\nu_0$).    On the other hand,  any $a_2 \in \{50, 100\}$ gives  $\mbox{ape}_\theta \leq 10\%$, for all combinations of $d_0 \in\{0.10, 0.25, 0.35, 0.45\}$ and   $\nu_0 \in \{1.1, 1.5, 2.5, 5.0\}$.   The simulation results for $a_2 = 50$ ($b_2  \approx 158.333$) are illustrated  in Figure \ref{table:gaussian}.

Figure \ref{table:gaussian} shows  the  sample mean  (solid circle) and the 95\% credibility interval  (solid line)  for the sample obtained from the  posterior distribution of  $\nu, d, \theta, \gamma$ and  $\omega$ (respectively, from top to bottom), for each combination of $d_0$ and  $\nu_0$.   The true parameter values $\nu_0, d_0, \theta_0, \gamma_0$ and  $\omega_0$ are represented in the corresponding row by the dashed line.    The graphs related  to  $\theta, \gamma$ and $\omega$  (respectively,  the third, fourth and fifth rows, from top to bottom) consider  the same scale for all $d_0 \in\{0.10, 0.25, 0.35, 0.45\}$.   Also,  for the parameters  $\theta, \gamma$ and $\omega$,  there is one graph for each $d_0$ and,  for each one of these graphs,   the true value of $\nu_0$  is indicated in the $x$-axis.

From  Figure \ref{table:gaussian} one observes  that,  for $\nu$ and $\omega$, the conclusion regarding the estimation bias and the credibility intervals are basically the same as in  Case 1 (see Table \ref{table:uniform}).   On the other hand,  under   C4.1 of  Case 4),   $\mbox{ape}_{\eta_i} < 10\%$, for  all $i\in \{1,\cdots, 5\}$  and any  combination of $d_0 \in\{0.10, 0.25, 0.35, 0.45\}$ and   $\nu_0 \in \{1.1, 1.5, 2.5, 5.0\}$ (compare  the parameter $\theta$ in Table \ref{table:uniform} and  Figure \ref{table:gaussian}).    As in  Case 1, the bias for $\theta$ is always positive when $\nu_0 = 1.1$ (for any $d_0$) and negative when $\nu_0 > 1.1$ and the bias for the  parameters $d$ and $\gamma$ does not seem to follow any pattern.    Under the current scenario,  $d_0, \theta_0$ and $\gamma_0$ are all contained in the respective credibility intervals, for any combination of  $d_0 \in\{0.10, 0.25, 0.35, 0.45\}$ and   $\nu_0 \in \{1.1, 1.5, 2.5, 5.0\}$.

When the true value of  $\gamma_0$ is not  used  to choose  $b_2$ (scenario  C4.2) similar results to the ones in Figure \ref{table:gaussian} are still obtained  for some combinations of $(a_2, b_2)$.   Not surprisingly, the pairs $(a_2, b_2)$ which lead to good estimates are such that  $a_2(a_2+b_2)^{-1}$ (the mean $\mu_B$ of the prior distribution) is close to $\gamma_0$.   For instance, when  $(a_2, b_2) \in \{(100,  300), (100,350)\}$ ($\mu_B$ is, respectively, equal to 0.25 and 0.22, while  $\gamma_0 = 0.24$)  it is obtained $\mbox{ape}_\gamma < 10\%$ for all combinations of $d_0 \in\{0.10, 0.25, 0.35, 0.45\}$ and   $\nu_0 \in \{1.1, 1.5, 2.5, 5.0\}$.   The pair $(a_2,b_2) = (100,350)$ provides slightly better results  than $(a_2,b_2) = (100,300)$ only when $\nu_0 = 5$.   When $(a_2,b_2) = (100,280)$ (so $\mu_B \approx 0.26$),  $\mbox{ape}_\gamma \leq 10\%$ in 16 out of 20 cases  (in the remaining 4 cases,   $\mbox{ape}_\gamma$ does not exceed  13.4\%).

On the other hand,  choosing $a_2$ and $b_2$ such that  $a_2(a_2+b_2)^{-1}$  is close to the true $\gamma_0$ does not necessarily lead to good estimates.   For instance, if  $(a_2, b_2) = (5, 15)$  then $\mu_B = 0.25$,  as it is when  $(a_2, b_2)  =  (100, 300)$,  but  $\mbox{ape}_\gamma > 10\%$ in  6 out of 20 cases.     Also, it is not evident that the more distant $a_2(a_2+b_2)^{-1}$  is from $\gamma_0$, the worst is the estimation.  For instance,   by letting   $(a_2, b_2) \in  \{(3, 15), (100, 440), (10, 30), (10, 40), (20, 80), (100, 400),   (100, 270), (5, 10)\}$ then, respectively,   $\mu_B \in \{0.167, 0.185, 0.200,0.200,0.200,0.200, 0.270, 0.330\}$  and  it is observed that $\mbox{ape}_\gamma > 10\%$ in 14, 20, 4, 10, 12,  16, 13 and 5 out of 20 cases, respectively.

\begin{table}[!htp]
\centering
\caption{Summary for the sample obtained from posterior distributions considering Beta priors for $2d,-\theta$ and $\gamma$:  mean $\bar \eta_i$, standard deviation $\mbox{sd}_{\eta_i}$ and the 95\% credibility interval $CI_{0.95}(\eta_i)$ for the parameter $\eta_i$  in $\boldsymbol {\eta} =  (\nu, d, \theta, \gamma, \omega)' :=  (\eta_1, \cdots, \eta_5)' $, for each  $i\in \{1,\cdots, 5\}$.  The true parameter values considered in this simulation are $d_0 \in\{0.10, 0.25, 0.35, 0.45\}$, $\nu_0 \in \{1.1, 1.5, 2.5, 5.0\}$, $\theta_0 = -0.15$, $\gamma_0 = -0.24$ and $\omega_0 = -5.4$.}\label{table:beta}
{\footnotesize
\begin{tabular}{c|ccccccccccc}
\toprule
\multicolumn{1}{c}{\multirow{2}{*}{$d_0$}}&
\multicolumn{1}{c}{\multirow{2}{*}{$\nu_0$}}&&
\multicolumn{1}{c}{$\bar{\nu}$ $(\mbox{sd}_{\nu})$} &&
\multicolumn{1}{c}{$\bar{d}$ $(\mbox{sd}_{d})$} &&
\multicolumn{1}{c}{$\bar{\theta}$ $(\mbox{sd}_{\theta})$} &&
\multicolumn{1}{c}{$\bar{\gamma}$ $(\mbox{sd}_{\gamma})$} &&
\multicolumn{1}{c}{$\bar{\omega}$ $(\mbox{sd}_{\omega})$} \\
\multicolumn{1}{c}{}&\multicolumn{1}{c}{}&&
\multicolumn{1}{c}{$CI_{0.95}(\nu)$} &&
\multicolumn{1}{c}{$CI_{0.95}(d)$} &&
\multicolumn{1}{c}{$CI_{0.95}(\theta)$} &&
\multicolumn{1}{c}{$CI_{0.95}(\gamma)$} &&
\multicolumn{1}{c}{$CI_{0.95}(\omega)$} \\

\cmidrule{1-2} \cmidrule{4-4} \cmidrule{6-6} \cmidrule{8-8}\cmidrule{10-10}\cmidrule{12-12}
\multirow{10}{*}{0.10}
 & \multirow{2}{*}{1.1} && 1.102 (0.045) &&0.100 (0.018) &&-0.144 (0.011) &&0.243 (0.026) &&-5.472 (0.063) \\
&&& [1.012; 1.198] && [0.069; 0.139] && [-0.170; -0.124] && [0.199; 0.294] && [-5.569; -5.347] \\
\cmidrule{2-2} \cmidrule{4-4} \cmidrule{6-6} \cmidrule{8-8}\cmidrule{10-10}\cmidrule{12-12}
 & \multirow{2}{*}{1.5} && 1.483 (0.063) &&0.102 (0.019) &&-0.154 (0.013) &&0.239 (0.024) &&-5.412 (0.047) \\
&&& [1.364; 1.607] && [0.069; 0.141] && [-0.178; -0.132] && [0.190; 0.285] && [-5.524; -5.333] \\
\cmidrule{2-2} \cmidrule{4-4} \cmidrule{6-6} \cmidrule{8-8}\cmidrule{10-10}\cmidrule{12-12}
 & \multirow{2}{*}{1.9} && 2.067 (0.108) &&0.101 (0.016) &&-0.160 (0.012) &&0.234 (0.027) &&-5.422 (0.034) \\
&&& [1.879; 2.309] && [0.070; 0.134] && [-0.186; -0.137] && [0.184; 0.291] && [-5.498; -5.360] \\
\cmidrule{2-2} \cmidrule{4-4} \cmidrule{6-6} \cmidrule{8-8}\cmidrule{10-10}\cmidrule{12-12}
 & \multirow{2}{*}{2.5} && 2.702 (0.143) &&0.106 (0.019) &&-0.160 (0.012) &&0.231 (0.023) &&-5.407 (0.022) \\
&&& [2.442; 2.977] && [0.074; 0.150] && [-0.186; -0.136] && [0.187; 0.277] && [-5.470; -5.375] \\
\cmidrule{2-2} \cmidrule{4-4} \cmidrule{6-6} \cmidrule{8-8}\cmidrule{10-10}\cmidrule{12-12}
 & \multirow{2}{*}{5.0} && 5.251 (0.391) &&0.099 (0.017) &&-0.158 (0.012) &&0.263 (0.024) &&-5.343 (0.036) \\
&&& [4.514; 6.080] && [0.070; 0.132] && [-0.186; -0.134] && [0.220; 0.312] &&  \textbf{[-5.390; -5.265]} \\
\midrule

\multirow{10}{*}{0.25}
 & \multirow{2}{*}{1.1} && 1.095 (0.044) &&0.255 (0.031) &&-0.145 (0.012) &&0.242 (0.028) &&-5.455 (0.057) \\
&&& [0.986; 1.190] && [0.197; 0.315] && [-0.173; -0.123] && [0.191; 0.295] && [-5.559; -5.348] \\
\cmidrule{2-2} \cmidrule{4-4} \cmidrule{6-6} \cmidrule{8-8}\cmidrule{10-10}\cmidrule{12-12}
 & \multirow{2}{*}{1.5} && 1.485 (0.066) &&0.252 (0.030) &&-0.154 (0.013) &&0.240 (0.022) &&-5.421 (0.046) \\
&&& [1.357; 1.628] && [0.193; 0.313] && [-0.180; -0.132] && [0.195; 0.281] && [-5.529; -5.324] \\
\cmidrule{2-2} \cmidrule{4-4} \cmidrule{6-6} \cmidrule{8-8}\cmidrule{10-10}\cmidrule{12-12}
 & \multirow{2}{*}{1.9} && 2.058 (0.113) &&0.260 (0.030) &&-0.159 (0.012) &&0.236 (0.028) &&-5.431 (0.048) \\
&&& [1.864; 2.316] && [0.194; 0.309] && [-0.185; -0.138] && [0.189; 0.293] && [-5.525; -5.340] \\
\cmidrule{2-2} \cmidrule{4-4} \cmidrule{6-6} \cmidrule{8-8}\cmidrule{10-10}\cmidrule{12-12}
 & \multirow{2}{*}{2.5} && 2.719 (0.146) &&0.271 (0.029) &&-0.163 (0.012) &&0.229 (0.021) &&-5.386 (0.038) \\
&&& [2.452; 3.008] && [0.216; 0.326] && [-0.188; -0.137] && [0.186; 0.275] && [-5.469; -5.304] \\
\cmidrule{2-2} \cmidrule{4-4} \cmidrule{6-6} \cmidrule{8-8}\cmidrule{10-10}\cmidrule{12-12}
 & \multirow{2}{*}{5.0} && 5.208 (0.326) &&0.244 (0.028) &&-0.159 (0.013) &&0.260 (0.023) &&-5.313 (0.034) \\
&&& [4.548; 5.871] && [0.190; 0.299] && [-0.185; -0.134] && [0.213; 0.309] &&  \textbf{[-5.381; -5.256]} \\
\midrule

\multirow{10}{*}{0.35}
 & \multirow{2}{*}{1.1} && 1.097 (0.038) &&0.355 (0.034) &&-0.145 (0.012) &&0.241 (0.027) &&-5.469 (0.080) \\
&&& [1.014; 1.175] && [0.283; 0.413] && [-0.169; -0.126] && [0.186; 0.298] && [-5.630; -5.289] \\
\cmidrule{2-2} \cmidrule{4-4} \cmidrule{6-6} \cmidrule{8-8}\cmidrule{10-10}\cmidrule{12-12}
 & \multirow{2}{*}{1.5} && 1.481 (0.064) &&0.349 (0.030) &&-0.154 (0.012) &&0.239 (0.024) &&-5.447 (0.077) \\
&&& [1.359; 1.628] && [0.285; 0.403] && [-0.178; -0.131] && [0.192; 0.285] && [-5.587; -5.310] \\
\cmidrule{2-2} \cmidrule{4-4} \cmidrule{6-6} \cmidrule{8-8}\cmidrule{10-10}\cmidrule{12-12}
 & \multirow{2}{*}{1.9} && 2.070 (0.102) &&0.370 (0.029) &&-0.160 (0.012) &&0.236 (0.027) &&-5.414 (0.082) \\
&&& [1.870; 2.288] && [0.306; 0.421] && [-0.183; -0.136] && [0.186; 0.293] && [-5.587; -5.237] \\
\cmidrule{2-2} \cmidrule{4-4} \cmidrule{6-6} \cmidrule{8-8}\cmidrule{10-10}\cmidrule{12-12}
 & \multirow{2}{*}{2.5} && 2.720 (0.147) &&0.375 (0.028) &&-0.163 (0.011) &&0.228 (0.023) &&-5.337 (0.063) \\
&&& [2.449; 3.009] && [0.321; 0.426] && [-0.185; -0.143] && [0.186; 0.274] && [-5.440; -5.221] \\
\cmidrule{2-2} \cmidrule{4-4} \cmidrule{6-6} \cmidrule{8-8}\cmidrule{10-10}\cmidrule{12-12}
 & \multirow{2}{*}{5.0} && 5.147 (0.346) &&0.344 (0.027) &&-0.159 (0.012) &&0.258 (0.023) &&-5.255 (0.048) \\
&&& [4.598; 5.965] && [0.287; 0.395] && [-0.185; -0.137] && [0.214; 0.305] &&  \textbf{[-5.321; -5.171]} \\
\midrule

\multirow{10}{*}{0.45}
 & \multirow{2}{*}{1.1} && 1.101 (0.042) &&0.454 (0.024) &&-0.145 (0.012) &&0.238 (0.026) &&-5.424 (0.128) \\
&&& [1.024; 1.191] && [0.402; 0.489] && [-0.169; -0.122] && [0.189; 0.286] && [-5.682; -5.160] \\
\cmidrule{2-2} \cmidrule{4-4} \cmidrule{6-6} \cmidrule{8-8}\cmidrule{10-10}\cmidrule{12-12}
 & \multirow{2}{*}{1.5} && 1.493 (0.073) &&0.450 (0.024) &&-0.154 (0.012) &&0.243 (0.024) &&-5.414 (0.132) \\
&&& [1.362; 1.645] && [0.395; 0.488] && [-0.177; -0.132] && [0.195; 0.291] && [-5.681; -5.139] \\
\cmidrule{2-2} \cmidrule{4-4} \cmidrule{6-6} \cmidrule{8-8}\cmidrule{10-10}\cmidrule{12-12}
 & \multirow{2}{*}{1.9} && 2.045 (0.109) &&0.464 (0.019) &&-0.158 (0.011) &&0.239 (0.026) &&-5.424 (0.126) \\
&&& [1.844; 2.241] && [0.419; 0.491] && [-0.181; -0.134] && [0.193; 0.291] && [-5.717; -5.216] \\
\cmidrule{2-2} \cmidrule{4-4} \cmidrule{6-6} \cmidrule{8-8}\cmidrule{10-10}\cmidrule{12-12}
 & \multirow{2}{*}{2.5} && 2.742 (0.142) &&0.466 (0.017) &&-0.164 (0.012) &&0.228 (0.022) &&-5.170 (0.084) \\
&&& \textbf{[2.507; 3.010]} && [0.431; 0.493] && [-0.191; -0.141] && [0.189; 0.275] &&  \textbf{[-5.308; -5.002]} \\
\cmidrule{2-2} \cmidrule{4-4} \cmidrule{6-6} \cmidrule{8-8}\cmidrule{10-10}\cmidrule{12-12}
 & \multirow{2}{*}{5.0} && 5.164 (0.346) &&0.447 (0.022) &&-0.158 (0.011) &&0.256 (0.022) &&-5.070 (0.076) \\
&&& [4.558; 5.883] && [0.396; 0.486] && [-0.182; -0.137] && [0.214; 0.300] &&  \textbf{[-5.227; -4.942]} \\

\bottomrule
\end{tabular}}
\begin{description}
\item[{\footnotesize {\bf Note:}}]{\footnotesize The bold-face font for the \textbf{credibility interval} indicates that the interval does not contain the true parameter value.}
\end{description}
\end{table}

\subsubsection*{Case 5:  Beta Priors for $2d$, $-\theta$ and $\gamma$.}

Analogously to all other cases, the  estimation of $\nu$, $\theta$, $\gamma$ and  $\omega$ is not significantly affected by the change in the prior for $d$.

Upon assuming $-\theta \sim \mbox{Beta}(a_1, b_1)$ and $\gamma \sim \mbox{Beta}(a_2, b_2)$,  with  the same  $a_1$,  $a_2$, $b_1$ and $b_2$ as in scenario  C4.1 of  Case 4, and letting $2d \sim  \mbox{Beta}(a_3, b_3)$, with  $b_3 = a_3(1-2d_0)(2d_0) ^{-1}$, for each  $d_0 \in\{0.10, 0.25, 0.35, 0.45\}$ (scenario  C5.1 of  Case 5),  the following is concluded.   If  $a_3 \in \{25, 50\}$ then $\mbox{ape}_d < 10\%$  for all $d_0 \in \{0.10, 0.25, 0.35, 0.45 \}$.   By  increasing or decreasing too much the $a_3$ values the estimation performance decays. For instance,  $a_3 \in \{0.10, 0.20, 2.00\}$ yields $\mbox{ape}_d > 10\%$ in 3, 1 and 7 cases, respectively.

 Table  \ref{table:beta} reports the simulation results for   $a_3 = 25$  and  $b_3 = a_3(1-2d_0)(2d_0) ^{-1}$, which gives $b_1 \in \{100.000, 25.000, 10.714, 2.778\}$, respectively,  for $d_0 \in\{0.10, 0.25, 0.35, 0.45\}$.   The conclusions on the results presented in this table are the same as in Figure \ref{table:gaussian}.   Although  the credibility intervals for $d$  are slightly wider in Table \ref{table:beta} than in Figure \ref{table:gaussian}, in both tables  $d_0 \in CI_{0.95}(d)$ for any combination of  $d_0 \in\{0.10, 0.25, 0.35, 0.45\}$ and   $\nu_0 \in \{1.1, 1.5, 2.5, 5.0\}$.

As in  Case 2, when  considering a two step estimator (scenario  C5.2),  no improvement is observed, when compared to   Case 1.   In fact, once again,   the  estimates obtained by letting $a_3 = 25$ and $b_3 = a_3(1-2\bar d)(2\bar d) ^{-1}$, where $\bar d$ is the estimate of $d$ obtained in  Case 1,  are very close to $\bar d$ itself.

\section{Conclusions}\label{seciton:conclusion}

The Bayesian inference approach for parameter estimation on FIEGARCH models was described and a Monte Carlo simulation study was conducted  to analyze  the performance of the method  under the presence of long-memory in volatility.  The samples  from FIEGARCH processes were obtained by considering the infinite sum representation for the logarithm of the volatility. A recurrence formula was used to obtain the coefficients for this representation.  The generalized error distribution, with different  tail-thickness parameters was considered   so both innovation processes with  lighter and heavier tails than the  Gaussian  distribution,   were covered.

Markov Chain Monte Carlo (MCMC) methods where used to obtain samples from the posterior distribution of the parameters. A sensitivity analysis was performed by considering the following steps.    First,   an improper prior for $\nu$ and uniform priors $d, \theta, \gamma$ and $\omega$ were selected. In this case,   only the basic set of information  usually available in practice was considered. Second,   non-uniform priors were selected for one or more parameters in $\{d, \theta, \gamma\}$.   A  Gaussian prior for $\phi^{-1}(d)$, with $\phi(\cdot)$ defined in \eqref{eq:phi}, combined with uniform or Beta priors for $\theta$ ($-\theta$ in the Beta case) and $\gamma$ was considered.  In the sequel,  a comparison was made by assuming  Beta priors for $2d$, $-\theta$ and $\gamma$.  The sensitivity analysis was completed by  integrating (or not) the knowledge on the true parameter values to select  the hyperparameter values.

An example was presented to illustrate the similarities or differences on the mean, standard deviation and credibility intervals estimated by considering a chain of size $N = 200801$,  a thinned chain  (thinning parameter 200 and burn-in size 1000) and a sample of size 1000 (obtained from the larger chain, after the burn-in of size 1000).   Given the ergodicity of the Markov chain, the posterior means for all three chains were very close.  The differences on the standard deviations and credibility intervals are not significant enough to justify the use of the entire or thinned chains.    Although the example only presents the case $d_0 = 0.25$, the same conclusions apply to $d_0 \in \{0.10,  0.35, 0.45\}$.

The simulation study showed that  if   the prior of one or more parameters  is changed,  the estimation of the other parameters is not significantly affected.    The parameters $\nu$ and $\omega$  are always well estimated, in terms of absolute percentage error,  regardless  priors considered for   $d, \theta$ and $\gamma$,  for any combination of $d_0 \in\{0.10, 0.25, 0.35, 0.45\}$ and   $\nu_0 \in \{1.1, 1.5, 2.5, 5.0\}$.    With a few exceptions,  the true parameter value $\nu_0$ was  contained in the 95\% credibility interval,  for any combination of  $\nu_0 \in \{1.1, 1.5, 1.9, 2.5, 5.0\}$ and $d_0 \in \{0.10, 0.25, 0.35, 0.45\}$ considered.  The true parameter value  $\omega_0$ was not contained in any credibility interval when $\nu_0 = 5$.

Regardless the prior considered,  the parameter $d$ is usually better estimated when  $d \in \{$0.35, 0.45$\}$.   The Gaussian prior for  $\phi^{-1}(d)$ only provided  better results (globally) when the knowledge on the true parameter value $d_0$ was used to set  $\mu_\phi = \phi^{-1}(d_0)$ and  $\sigma_\phi$ was set to some value  smaller or equal than 1.  In particular,  only when $b = 0.15$ the absolute percentage error of estimation (ape) became smaller than 10\% for all $d_0 \in\{0.10, 0.25, 0.35, 0.45\}$.   Although  the credibility intervals for $d$  are slightly wider when a Beta prior is considered,  the use of the Beta prior for $2d$ neither improves nor degrades the estimation performance, compared to the Gaussian prior for $ \phi^{-1}(d)$.

The absolute percentage error of  estimation  for  $\theta$ ($\mbox{ape}_\theta$) only became smaller than 10\% when the Beta prior was considered and the true value of the parameter was used to select the hyperparameter.  When   $\theta_0$ was assumed unknown the $\mbox{ape}_\theta$ was always between 10\% and  38.1\%.  The parameter $\gamma$ is always better estimated than $\theta$, for any priors considered.    Similar to $d$ and $\theta$, the best performance is obtained when  the true parameter value is used to select the hyperparameters.   On the other hand,  $\gamma$ is the only parameter for which there are hyperparameter values  that do not yield $\mu_B = \gamma_0$ ($\mu_B$ is the mean of the prior distribution and $\gamma_0$ is the true parameter value)  while still providing good estimates.

 \section*{Acknowledgments}

 T.S. Prass was supported by CNPq-Brazil.  S.R.C. Lopes was partially supported by CNPq-Brazil, by CAPES-Brazil, by INCT em
 \emph{Ma\-te\-m\'a\-ti\-ca} and by Pronex {\it Probabilidade e Processos Estoc\'asticos} - E-26/170.008/2008 -APQ1.
 The authors are grateful to the (Brazilian) National Center of Super Computing (CESUP-UFRGS) for the computational resources.


\begin{thebibliography}{10}



\bibitem[{{Baillie et al.}(1996)}]{BEA1996}  Baillie, R.;  T. Bollerslev and H.O. Mikkelsen (1996). ``Fractionally Integrated Generalized Autoregressive Conditional Heteroskedasticity''. \emph{Journal of Econometrics}, vol. 74,  3-30.

\bibitem[{{Bollerslev}(1986)}]{B1986}  Bollerslev, T. (1986). ``Generalized Autoregressive Conditional Heteroskedasticity''. \emph{Journal of Econometrics}, vol. 31, 307-327.


\bibitem[{{Bollerslev and Mikkelsen}(1996)}]{BM1996} Bollerslev, T. and H.O. Mikkelsen (1996).  ``Modeling and Pricing Long  Memory in Stock Market Volatility". \emph{Journal of Econometrics}, vol.  73,  151-184.

\bibitem[{{Breidt et al.}(1998)}]{BEA1998}  Breidt, F.;  N. Crato and  P.J.F. de Lima (1998).  ``On the Detection and Estimation of Long Memory in Stochastic Volatility''. \emph{Journal of Econometrics},  vol.  83, 325-348.


 \bibitem[{{Brockwell and Davis}(1991)}]{BD1991} Brockwell,  P.J.  and R.A. Davis (1991). \emph{Time Series: Theory
      and Methods}. Second Edition. New York: Springer-Verlag.

\bibitem[{{Casela and George}(1992)}]{CG1992} Casella,  G. and E.I. George (1992). ``Explaining the Gibbs Sampler", \emph{The American Statistician},  vol  46(3), 167-174.


\bibitem[{{Chib and Greenberg}(1995)}]{CG1995}  Chib, S. and E. Greenberg (1995). ``Understanding the Metropolis-Hastings Algorithm", \emph{The American Statistician}, vol. 49(4), 327-335.


\bibitem[{{Engle}(1982)}]{E1982}  Engle, R.F. (1982). ``Autoregressive Conditional Heteroskedasticity with Estimates of Variance of U.K. Inflation". \emph{Econometrica}, vol. 50, 987-1008.


\bibitem[{{Gelfand and Smith}(1990)}]{GS1990} Gelfand, A.E. and A.F.M Smith  (1990). ``Sampling-Based Approaches to Calculating Marginal Densities".  \emph{Journal of the American Statistical Association}, vol. 85(410),  398-409.

\bibitem[{{Geman and Geman}(1984)}]{GG1984} Geman, S. and D. Geman (1984).  ``Stochastic Relaxation, Gibbs Distributions and the Bayesian Restoration of Images".  \emph{IEEE Transactions on Pattern Analysis and Machine Intelligence}, vol. 12, 609-628.

\bibitem[{{Gelman and Rubin}(1992)}]{GR1992} Gelman, A. and D. Rubin (1992).  ``Inference from Iterative Simulation Using Multiple Sequences". \emph{Statistical Science}, vol. 7, 457-511.

\bibitem[{{Geyer}(1992)}]{G1992} Geyer, C. J. (1992).  ``Practical Markov chain Monte Carlo". \emph{Statistical Science}, vol. 7, 473-511.

\bibitem[{{Hastings}(1970)}]{H1970} Hastings, W.K. (1970). ``Monte Carlo Sampling Methods Using Markov Chains and Their Applications".  \emph{Biometrika}, vol. 57(1), 97-109.

\bibitem[{{Lopes and Prass}(2013)}]{LP2013} Lopes, S.R.C  and  T.S. Prass (2013). ``Theoretical Results on  Fractionally Integrated Exponential   Generalized Autoregressive Conditional Heteroskedastic  Processes." Working Paper.


  \bibitem[{{MacEachern and Berliner}(1994)}]{MB1994} MacEachern, S.N. and  L.M. Berliner (1994).  ``Subsampling the Gibbs Sampler".  \emph{The American  Statistician}, vol. 48, 188-190.

\bibitem[{{Metropolis et al.}(1953)}]{M1953}  Metropolis, N.; A.W. Rosenbluth; M.N. Rosenbluth; A.H. Teller and  E. Teller  (1953). ``Equations of State Calculations by Fast Computing Machines". \emph{Journal of Chemical Physics},  vol.  21(6),  1087-1092.

\bibitem[{{Meyer and Yu}(2000)}]{MY2000} Meyer, R. and T. Yu (2000). ``Bugs for a Bayesian Analysis of Stochastic Volatility Models". \emph{Econometrics Journal}, vol. 3, 198-215.

\bibitem[{{Nelson}(1991)}]{NE1991}  Nelson, D.B. (1991). ``Conditional Heteroskedasticity in Asset Returns: A
    New Approach''. \emph{Econometrica}, vol.  59,  347-370.

\bibitem[{{Ruiz and Veiga}(2008)}]{RV2008} Ruiz E. and  H. Veiga  (2008). ``Modelling long-memory volatilities with leverage effect: A-LMSV versus FIEGARCH''. \emph{Computational Statistics and Data Analysis}, vol. 52(6), 2846-2862.

\bibitem[{{Smith and Roberts}(1993)}]{SR1993} Smith, A.F.M. and G.O. Roberts  (1993). ``Bayesian computation via the Gibbs sampler and related Markov chain Monte Carlo methods". \emph{Journal of the Royal Statistical Society}, Ser. B, vol. 55, 3-23.




 \end{thebibliography}
\end{document}